\documentclass[11pt]{article}

\usepackage{amssymb}
\usepackage[mathscr]{euscript}

\newtheorem{rema}{Remarque}
\newtheorem{defi}{D\'efinition}
\newtheorem{lemm}{Lemme}
\newtheorem{theo}{Th\'eor\`eme}

\newtheorem{exem}{Exemple}

\newcommand{\R}[1][]{\ensuremath{{\mathbb{R}^{#1}} }}
\newcommand{\C}[1][]{\ensuremath{{\mathbb{C}^{#1}} }}

\def\<{\langle} \def\>{\rangle}

\newcommand{\1}{\mathrm{1 \hspace{-0.25em} l}}

\title{Probl\`emes variationnels invariants par transformation conforme en
dimension 2}
\author{Fr\'ed\'eric H\'elein,}
\date{Le 29 janvier 2001}

\begin{document}
\maketitle
\section{Introduction}
Soit $u$ une fonction d\'efinie sur un domaine
$\Omega$ inclus dans $\Bbb{C}\simeq \Bbb{R}^2$ et \`a valeurs dans \R. La fonctionnelle de Dirichlet

$${\cal E}[u]:= \int _{\Omega}\left( {\partial u\over \partial x}^2
+ {\partial u\over \partial y}^2\right) dxdy = \int _{\Omega}|\nabla
u|^2dxdy$$
est un exemple bien connu de probl\`eme variationnel invariant par
transformation
conforme.
Ses points critiques sont les fonctions harmoniques, solutions de l'\'equation
$\Delta u =0$ et sont \'etroitement reli\'es aux fonctions holomorphes,
puisque
$u$ est harmonique si et seulement si
$x+iy\longmapsto \partial _xu -i\partial _yu$ est holomorphe. L'\'etude des
surfaces minimales dans l'espace de dimension 3 a conduit assez
naturellement les math\'ematiciens \`a consid\'erer la m\^eme fonctionnelle
pour des applications d'un domaine $\Omega \subset \C$ dans \R[3]. En effet,
pour toute application $u:\Omega \longrightarrow \R[3]$, on a toujours
${\cal E}[u]\geq {\cal A}[u]$, o\`u ${\cal A}$ est la fonctionnelle aire
(rendue extr\'emale pour les surfaces minimales), d\'efinie par

$${\cal A}[u] := \int _{\Omega}\left| {\partial u\over \partial x}\times
{\partial u\over \partial y}\right| dxdy.$$
Et de plus, ${\cal A}[u]={\cal E}[u]$ si $u$ est conforme. Or la
fonctionnelle ${\cal E}$ poss\`ede
des propri\'et\'es de compacit\'e bien meilleures que ${\cal A}$. Cela a
permis
\`a J. Douglas et T. Rad\'o en 1930 de r\'esoudre le probl\`eme de Plateau. Le
succ\`es de cette approche a \'et\'e confirm\'e par les travaux de C.B.
Morrey:
on peut construire des surfaces minimales dans des vari\'et\'es riemanniennes
$({\cal N},g)$ en cherchant les points critiques - appel\'es {\em applications
harmoniques} - de

$${\cal E}[u]:= \int _{\Omega}g_{ij}(u)\left( {\partial u^i\over \partial x}
{\partial u^j\over \partial x} + {\partial u^i\over \partial y}
{\partial u^j\over \partial y}\right) dxdy,$$
qui g\'en\'eralise la fonctionnelle de Dirichlet classique donn\'ee plus haut.
On peut \'egalement modifier ${\cal E}[u]$ en y ajoutant un terme du type
$\int _{\Omega}u^{\star}\omega$, o\`u $\omega$ est une deux-forme
d\'efinie sur ${\cal N}$ et $u^{\star}\omega$ son image inverse par $u$.
Si ${\cal N}$ est de dimension 3, on obtient ainsi une fonctionnelle dont
la th\'eorie variationnelle produit des surfaces \`a courbure moyenne
prescrite (\'egale au rapport entre $d\omega$ et la forme de
de volume riemannien sur ${\cal N}$), cf [H\'e].

Toutes ces actions sont des exemples de fonctionnelles invariantes par
le groupe des transformations conformes de \R[2]. Ce sont m\^eme les seules
possibles, lorsque l'on se restreint \`a des fonctionelles du type
$\int _{\Omega}L(u,\nabla u)dxdy$, en supposant que $\xi \mapsto L(u,\xi )$
est quadratique.

L'int\'er\^et de ce type de fonctionnelle est aujourd'hui bien \'etabli
en physique math\'ematique, puisque la th\'eorie des cordes
et des supercordes est b\^atie sur une quantification de ${\cal E}$ et de
ses g\'en\'eralisations supersym\'etriques. En particulier, comme il a \'et\'e montr\'e A. M. Polyakov, l'\'energie ${\cal E}$
s'av\`ere
\^etre plus appropri\'ee que l'aire ${\cal A}$ pour le calcul
d'int\'egrales fonctionnelles.

Notre but ici, est de d\'ecrire une classe g\'en\'erale de probl\`emes
variationnels en dimension
2 qui sont invariants par le groupe des transformations conformes de \C .

Nous verrons en particulier qu'une g\'eom\'etrie,
que nous appelons $\Bbb{C}$-{\em finslerienne}, similaire \`a la g\'eom\'etrie
finslerienne est associ\'ee de fa\c{c}on naturelle
\`a ces probl\`emes. Essentiellement une m\'etrique $\Bbb{C}$-finslerienne est la donn\'ee d'une
application $F:T{\cal N}^{\Bbb{C}}\longrightarrow [0,\infty [$ homog\`ene de degr\'e deux sur chaque
fibre, c'est \`a dire telle que $F(y,\lambda z) = |\lambda|^2F(y,z)$, $\forall y\in {\cal N}$,
$\forall z\in T_y{\cal N}^{\Bbb{C}}$ et $\forall \lambda \in \Bbb{C}$.
A travers une analyse succincte de cette g\'eom\'etrie, nous verrons qu'apparemment elle partage
nombre de propri\'et\'es avec la g\'eom\'etrie finslerienne classique.

Dans une deuxi\`eme partie, nous nous int\'eressons aux formulations hamiltoniennes
pour les probl\`emes invariants par transformations conformes. Rappelons qu'en calcul
des variations \`a plusieurs variables, les possibilit\'es sont multiples.
Nous explorerons d'abord bri\`evement le formalisme de De Donder-Weyl, puis nous nous
int\'eresserons au formalisme de Carath\'eodory. Le lecteur un peu sp\'ecialiste de cette
th\'eorie constatera que j'ai introduit un param\`etre suppl\'ementaire, not\'e $w$
(qui doit \^etre remplac\'e par 0 si on veut comparer ce qui est \'ecrit ici avec, par
exemple, l'expos\'e dans [Ru]). Cela a \'et\'e
motiv\'e par le fait que l'analogue de la transformation de Legendre pour la th\'eorie
de Carath\'eodory est en g\'en\'eral mal d\'efini, sans ce degr\'e de libert\'e suppl\'ementaire.
Pour des d\'eveloppements suppl\'ementaires concernant ces diff\'erents types de formalisme
hamiltonien, voir [HK].\\

\noindent {\bf Remerciements} Je tiens \`a remercier Joseph Kouneiher pour
les discussions que j'ai eu avec lui sur ce sujet et ses encouragements.

\section{Caract\'erisation}
\subsection{Etude locale}
Nous cherchons, parmi tous les lagrangiens $L(t,u,du)$ continument
diff\'erentiables, ceux qui sont tels que
la fonctionnelle d\'efinie par

$${\cal L}[u] = \int _{\Omega} L(t, u, du )dt^1dt^2$$
est invariante par le groupe des transformations conformes de \C $\simeq$
\R[2].
Au pr\'ealable, il faut rappeler le sens que nous donnons \`a cette notion
d'invariance.
Consid\'erons une famille de diff\'eomorphismes locaux $\Psi _s$ de
$\Bbb{R}^2$, \`a un 
param\`etre $s$, qui forme un 
groupe pour la composition. $\Psi _s$ est le flot d'un champ de vecteurs $X$ 
d\'efini sur un ouvert de
$\R ^2$ contenant l'adh\'erence de $\Omega$. Cela entra{\^\i}ne en
particulier que pour $s$ 
proche de 0, on a

\begin{equation}\label{1.1}\Psi _s(t)=t+sX(t)+o(s).\end{equation}
L'image par $\Psi _s$ de $\Omega$ est un ouvert $\Omega _s$, diff\'erent de 
$\Omega$ en g\'en\'eral.
Soit maintenant une application $u$ de $\Omega$ vers \R[n]: nous dirons
qu'elle est transform\'ee en $u_s$ si le graphe de $u_s$ est 
l'image du graphe de $u$
par la transformation $(t,y)\longmapsto (\Psi_s(t),y)$ agissant sur
$\R[2]\times \R[n]$. Donc le 
domaine de d\'efinition
de $u_s$ sera $\Omega _s=\Psi _s(\Omega )$, et $u_s$ satisfait \`a

\begin{equation}\label{1.2}u_s\circ \Psi _s =u,\forall s.\end{equation}
Nous dirons que la fonctionnelle ${\cal L}$ est invariante 
par $X$ si et seulement
si pour tout sous-domaine $\omega \subset \Omega$,

$${\cal L}_{\Psi _s(\omega )}[u_s]={\cal L}_{\omega }[u].$$
Une fa\c{c}on d'\'ecrire cette relation est de faire le changement de
variable 
$t=\Psi _s(\tau )$, pour
$\tau \in \omega$ dans l'int\'egrale de gauche. Cela donne

$$\int _{\omega}L\left( \Psi _s(\tau ),u_s(\Psi _s(\tau )),
du_s(\Psi _s(\tau ))\right) det(d\Psi _s(\tau ))d\tau
={\cal L}_{\omega }[u].$$
Or, en d\'erivant la relation (\ref{1.2}), on obtient:

$$du_s(\Psi _s(\tau )).d\Psi _s(\tau )=du(\tau ),$$
d'o\`u

$$du_s(\Psi _s(\tau ))=du(\tau ).d\Psi _s(\tau )^{-1}.$$
Donc, en utilisant cette relation et (\ref{1.2}), on obtient

\begin{equation}\label{1.3}
\int _{\omega}L\left( \Psi _s(\tau ),u(\tau ), du(\tau ).
d\Psi_s(\tau)^{-1}\right)
\hbox{det}(d\Psi _s(\tau ))d\tau
={\cal L}_{\omega }[u].
\end{equation}
Nous pouvons d\'eduire une version infinit\'esimale de cette relation, en 
supposant que $s$ est petit
et en d\'eveloppant au premier ordre:

$${\cal L}_{\omega }[u]=\int _{\omega}L\left( t+sX(t),u(t),du(t).
(\1 -sdX(t))\right)
\hbox{det}(\1 +sdX(t))dt +o(s).$$
Et comme cette relation doit \^etre valable pour tout $\omega$, 
n\'ecessairement,
$\forall (t,y,z)\in \Omega \times \R ^n\times M(\R ^2,\R ^n)$,

\begin{equation}\label{1.6}
L\left( t+sX(t),y,z.(\1 -sdX(t))\right) (1 +s\hbox{div}X(t))=L(t, y, z)+o(s).
\end{equation}
De fa\c{c}on \'equivalente:
$\forall (t,y,z)\in \Omega \times \R ^n\times M(\R ^2,\R ^n)$,

\begin{equation}\label{1.4}
X^{\alpha}(t){\partial L\over \partial t^{\alpha}}(t,y,z)
- {\partial L\over \partial z^i_{\alpha}}(t,y,z) z^i_{\beta}
{\partial X^{\beta}\over \partial t^{\alpha}}(t)
+ L(t,y,z){\partial X^{\alpha}\over \partial t^{\alpha}}(t) = 0.
\end{equation}

La question est de trouver les conditions sur $L$ pour que cette relation soit
vraie pour tout groupe \`a un param\`etre d'applications conformes $\Psi_s$.
En testant (\ref{1.4}) avec, comme groupe de d\'eformations les
translations de \R[2], engendr\'ees par les champs de vecteur
constants, on obtient ${\partial L\over \partial t^{\alpha}}=0$ partout,
\`a savoir que $L$ ne d\'epend pas de $t$. Ainsi (\ref{1.4}) se simplifie en

\begin{equation}\label{1.5}
\left(  {\partial L\over \partial z^i_{\alpha}}(y,z) z^i_{\beta}
-  L(y,z)\delta ^{\alpha}_{\beta} \right)
{\partial X^{\beta}\over \partial t^{\alpha}}(t) = 0.
\end{equation}

De mani\`ere g\'en\'erale, chaque $\Psi_s$ satisfait les \'equations
de Cauchy-Riemann $\partial _{t^1}\Psi_s^1 - \partial _{t^2}\Psi_s^2 =
\partial _{t^2}\Psi_s^1 + \partial _{t^1}\Psi_s^2 =0$, donc, d'apr\`es
(\ref{1.1}), $X$ est aussi holomorphe, \`a savoir,

$$\partial _{t^1}X^1 - \partial _{t^2}X^2 =
\partial _{t^2}X^1 + \partial _{t^1}X^2 =0.$$
Et $L$ satisfait (\ref{1.5}) pour tout champ de vecteur holomorphe
si et seulement si

\begin{equation}\label{1.5b}
\left\{ \begin{array}{ccl}
\displaystyle {\partial L\over \partial z^i_1}(y,z) z^i_1
+ {\partial L\over \partial z^i_2}(y,z) z^i_2 & = & 2L(y,z)\\
\displaystyle {\partial L\over \partial z^i_1}(y,z) z^i_2
- {\partial L\over \partial z^i_2}(y,z) z^i_1 & = & 0.
\end{array} \right.
\end{equation}

Nous pouvons exprimer la condition (\ref{1.5b}) de deux fa\c{c}ons
diff\'erentes. Premi\`erement,
nous d\'efinissons le {\em tenseur hamiltonien}

$$H^{\alpha}_{\beta}(t)=\sum _{i=1}^n{\partial u^i\over \partial t^{\beta}}
{\partial L\over \partial z^i_{\alpha}}(t, u(t),du(t))-\delta 
^{\alpha}_{\beta}L(t, u(t),du(t)),$$
g\'en\'eralisant l'hamiltonien ou l'\'energie totale associ\'ee \`a un
probl\`eme variationnel de
dimension 1. Et il est clair que (\ref{1.5b}) \'equivaut aux relations
$H^1_1+H^2_2 = H^1_2-H^2_1 = 0$,
signifiant que le tenseur hamiltonien est sym\'etrique et \`a trace nulle.

Deuxi\`emement, en identifiant l'ensemble des variables
$\{z^i_{\alpha}/i=1,...,n; \alpha=1,2\}$ avec
$\Bbb{C}^n$, on peut d\'efinir le lagrangien comme une fonction $F$ des
variables $(y,z)\in \Bbb{R}^n\times \Bbb{C}^n$
en notant

$$L(u,du) = F(u, {\partial u\over \partial t^1} + i{\partial u\over
\partial t^2}).$$
Fixons $y\in \Bbb{R}^n$, $z = z_1+iz_2\in \Bbb{C}^n$ et diff\'erentions la
fonction 

$$\begin{array}{ccl}
\Bbb{C} & \longrightarrow & \Bbb{R}\\
\lambda & \longmapsto & F(y,\lambda z)
\end{array}$$
par rapport \`a $\lambda=a+ib$. En utilisant (\ref{1.5b}), il vient

$$\begin{array}{ccl}
dF(y,\lambda z) & = & \displaystyle
\left( {\partial L\over \partial z^i_1}(y,\lambda z)z^i_1 + {\partial
L\over \partial z^i_2}(y,\lambda z)z^i_2\right) da
+ \left( {\partial L\over \partial z^i_2}(y,\lambda z)z^i_1 - {\partial
L\over \partial z^i_1}(y,\lambda z)z^i_2\right) db\\
& = & \displaystyle 2L(y,\lambda z){ada+bdb\over a^2+b^2} = F(y,\lambda
z){d|\lambda |^2\over |\lambda |^2}.
\end{array}$$
d'o\`u ${d\over d\lambda}\left( |\lambda |^{-2}F(y,\lambda z)\right) =0$,
c'est \`a dire $|\lambda |^{-2}F(y,\lambda z)$
ne d\'epend pas de $\lambda$. Nous en d\'eduisons le r\'esultat suivant.

\begin{theo}
L'action ${\cal L}[u]:= \int L(t,u(t),du(t))dt^1dt^2$ est invariante par
transformations conformes si et seulement si
$L(t,y^j,z^j_{\alpha}) = F(y^j,z^j_1+iz^j_2)$, o\`u $F:\Bbb{R}^n\times
\Bbb{C}^n\longrightarrow \Bbb{R}$
satisfait $\forall \lambda \in \Bbb{C}$, $\forall y\in \Bbb{R}^n$, $\forall
z\in \Bbb{C}^n$

\begin{equation}\label{1.6}
F(y,\lambda z) = |\lambda |^2F(y, z).
\end{equation}
\end{theo}

\subsection{Lois de conservation}
L'action ${\cal L}$ \'etant en particulier invariante par translations, le
th\'eor\`eme de Noether permet
de d\'eduire que le tenseur hamiltonien est \`a divergence nulle:
${\partial H^{\alpha}_{\beta}\over \partial x^{\alpha}}=0$, pour $\beta =1,
2$. Comme de plus
$H^{\alpha}_{\beta}$ est sym\'etrique \`a trace nulle, nous pouvons
reformuler cette loi de conservation en introduisant la
{\em diff\'erentielle de Hopf g\'en\'eralis\'ee} ${\cal Q}:=f(dz)^2$, avec

$$\begin{array}{ccl}
f & = & \left( (H^1_1-H^2_2) - i(H^1_2+H^2_1)\right) \\
& = & \displaystyle 
\left( {\partial L\over \partial z_1^j}(u,du){\partial u^j\over \partial
t^1} -
{\partial L\over \partial z_2^j}(u,du){\partial u^j\over \partial t^2}\right) 
-i \left( {\partial L\over \partial z_1^j}(u,du){\partial u^j\over \partial
t^2} +
{\partial L\over \partial z_2^j}(u,du){\partial u^j\over \partial
t^1}\right) ,
\end{array}$$
et en \'ecrivant que ${\cal Q}$ est holomorphe:

$${\partial f\over \partial \overline{z}} = 0.$$

\subsection{Un point de vue g\'eom\'etrique}

Nous pouvons ais\'ement g\'en\'eraliser ce qui pr\'ec\`ede \`a des
applications \`a valeurs dans une vari\'et\'e
${\cal N}$. L'action $\int L(t,u,\partial _1u,\partial _2u)dt^1dt^2$ est
invariante par transformation conforme
si et seulement si $L(t,u,\partial _1u,\partial _2u) = F(u,\partial
_1u+i\partial _2u)$, o\`u $F$ est
cette fois une application d\'efinie sur le fibr\'e tangent complexifi\'e

$$T{\cal N}^{\Bbb{C}} = T{\cal N}\otimes \Bbb{C} =
\{ (y,z)/y\in {\cal N},z\in T_y{\cal N}\otimes \Bbb{C} \},$$
telle que $\forall \lambda \in \Bbb{C}$, $F(y,\lambda z) = |\lambda
|^2F(y,z)$.

En imitant la d\'efinition d'une vari\'et\'e de Finsler, nous sommes
naturellement conduits \`a introduire ce qui
suit:

\begin{defi}
Une {\em pseudo-vari\'et\'e $\Bbb{C}$-finslerienne} est une vari\'et\'e
diff\'erentielle ${\cal N}$ munie d'une
application $F:T{\cal N}^{\Bbb{C}}\longrightarrow \Bbb{R}$ satisfaisant la
condition

$$\forall \lambda \in \Bbb{C},\ F(y,\lambda z) = |\lambda |^2F(y, z).$$
Une {\em vari\'et\'e $\Bbb{C}$-finslerienne} est une pseudo-vari\'et\'e
$\Bbb{C}$-finslerienne telle que
$F$ satisfait les conditions suivantes: $F$ est deux fois diff\'erentiable sur
$T{\cal N}^{\Bbb{C}}\setminus ({\cal N}\times \{0\})$
et il existe une constante $c>0$, telle que $\forall (y,z)\in T{\cal
N}^{\Bbb{C}}$,
$\forall v=v_1+iv_2\in T_y{\cal N}^{\Bbb{C}}$,

\begin{equation}\label{1.7}
{\partial ^2F\over \partial z^i_{\alpha}\partial
z^j_{\beta}}(y,z)v^i_{\alpha}v^j_{\beta} \geq c|v|^2.
\end{equation}
\end{defi}
Les cons\'equences de la ``condition d'ellipticit\'e'' (\ref{1.7}) seront
donn\'ees plus loin.
Si ${\cal N}$ est une vari\'et\'e $\Bbb{C}$-finslerienne, pour toute
surface de Riemann $\Sigma$ et pour
toute application $u:\Sigma \longrightarrow {\cal N}$, nous pouvons
d\'efinir son ``\'energie''

$${\cal E} [u]= \int _{\Sigma} F(u,2\overline{\partial} u)d\sigma ,$$
o\`u, dans toute coordonn\'ee locale holomorphe $t=t^1+it^2$ sur $\Sigma$,
$2\overline{\partial} u = \partial _1u+i\partial _2u$ et $d\sigma
=dt^1\wedge dt^2$.

\begin{exem}
Toute vari\'et\'e riemannienne ${\cal N}$ est une vari\'et\'e
$\Bbb{C}$-finslerienne. Si $g$ est le tenseur m\'etrique,
la fonction $F$ correspondante est juste $F(y,z) = {1\over
2}g_{ij}(y)(z^i_1z^j_1+ z^i_2z^j_2)$. La m\^eme
construction \`a partir d'une vari\'et\'e pseudo-riemannienne (telle que
$g_{ij}$ soit une m\'etrique de Minkowski)
donne une pseudo vari\'et\'e $\Bbb{C}$-finslerienne. Enfin si, en plus du
tenseur m\'etrique $g_{ij}$, on se donne
un tenseur antisym\'etrique $\omega _{ij}$ sur ${\cal N}$, on obtient une
vari\'et\'e $\Bbb{C}$-finslerienne
avec $F(y,z) = {1\over 2}\left( g_{ij}(y)(z^i_1z^j_1+ z^i_2z^j_2) + \omega
_{ij}(y)(z^i_1z^j_2- z^i_2z^j_1)\right)$.
\end{exem}

\section{La transform\'ee de Legendre}
Dans tout ce qui suit, nous supposerons que $F$ est deux fois
diff\'erentiable sur
$T{\cal N}^{\Bbb{C}}\setminus ({\cal N}\times \{0\})$ et nous noterons,
pour $\alpha ,\beta =1,2$ et $i,j=1,...,n$,

$$G^{\beta \alpha}_{ji}(y,z) = G^{\alpha \beta}_{ij}(y,z):= {\partial
^2F\over \partial z^i_{\alpha}\partial z^j_{\beta}}(y,z).$$

\subsection{Hypoth\`ese de Legendre}
Nous aurons besoin dans la suite, de supposer\\

\noindent {\bf Hypoth\`ese de Legendre globale} {\em Pour tout $y\in {\cal
N}$ et pour
tout $p=p^1+ip^2\in T_y^{\star}{\cal N}^{\Bbb{C}}$,
il existe un unique $z=z_1+iz_2\in T_y{\cal N}^{\Bbb{C}}$, tel que}

\begin{equation}\label{2.1}
{\partial F\over \partial z_1^j}(y,z) + i{\partial F\over \partial
z_2^j}(y,z) = p^1_j+ip^2_j .
\end{equation}
Cette propri\'et\'e peut \^etre montr\'ee localement, en utilisant le
th\'eor\`eme d'inversion locale, si l'on suppose ce qui suit:\\

\noindent {\bf Hypoth\`ese de Legendre locale} {\em Pour tout
$a=a^1+ia^2\in T_y^{\star}{\cal N}^{\Bbb{C}}$,
il existe un unique $v=v_1+iv_2\in T_y{\cal N}^{\Bbb{C}}$, tel que}

\begin{equation}\label{2.2}
G^{\alpha \beta}_{ij}(y,z)v^j_{\beta} = a^{\alpha}_i.
\end{equation} 

Mais dans le cas des vari\'et\'es $\Bbb{C}$-finsleriennes, la condition
d'ellipticit\'e (\ref{1.7}) entra\^\i ne l'hypoth\`ese
de Legendre globale (\ref{2.1}). En effet, cela implique que, pour tout
$y\in {\cal N}$ fix\'e, la fonction
$z\longmapsto F(y,z)$ est strictement convexe, puisqu'alors, pour tout
$v,w\in T_y{\cal N}^{\Bbb{C}}$ la
d\'eriv\'ee seconde de $l(s) = F(y,(1-s)v+xw)$ est $l''(s)\geq c|w-v|^2$.
Donc, pour tout $p\in T_y^{\star}{\cal N}^{\Bbb{C}}$,
une unique solution \`a (\ref{2.1}) est obtenue en minimisant $z\longmapsto
F(y,z)-p^{\alpha}_iz^i_{\alpha}$.
De plus, par homog\'en\'eit\'e, $F(y,0) = {\partial F\over \partial
z^j_{\alpha}}(y,0) = 0$ et donc la convexit\'e
de $F$ entra\^\i ne $F(y,z)\geq 0$, $\forall y,z$.

\subsection{Impulsions g\'en\'eralis\'ees}
Nous notons $2{\partial F\over \partial \overline{z}}:T{\cal
N}^{\Bbb{C}}\longrightarrow T^{\star}{\cal N}^{\Bbb{C}}$ l'application
d\'efinie par
$$2{\partial F\over \partial \overline{z}^j}(y,z) := {\partial F\over
\partial z_1^j}(y,z) + i{\partial F\over \partial z_2^j}(y,z).$$
Ses parties r\'eelles et imaginaires sont les impulsions g\'en\'eralis\'ees.
Nous faisons ici l'hypoth\`ese de Legendre globale. Alors nous d\'efinissons
$\psi_1(y,p^1,p^2) := z_1$ et $\psi_2(y,p^1,p^2) := z_2$ comme \'etant les
solutions de  l'\'equation (\ref{2.1}).
Nous notons $\Psi : T^{\star}{\cal N}^{\Bbb{C}}\longrightarrow T{\cal
N}^{\Bbb{C}}$ l'application
d\'efinie par

$$\Psi (y, p^1+ip^2) := \psi_1(y,p^1,p^2) + i\psi_2(y,p^1,p^2).$$
Donc $\Psi$ est l'application inverse de $2{\partial F\over \partial
\overline{z}}$.
\begin{lemm}
Les applications ${\partial F\over \partial \overline{z}}$ et $\Psi$
satisfont les relations suivantes:\\
$\forall (y,z)\in T{\cal N}^{\Bbb{C}}$, $\forall \lambda \in \Bbb{C}$,
\begin{equation}\label{2.3}
{\partial F\over \partial \overline{z}} (y,\lambda z) = \lambda {\partial
F\over \partial \overline{z}} (y,z),
\end{equation}
$\forall (y,p)\in T^{\star}{\cal N}^{\Bbb{C}}$, $\forall \lambda \in \Bbb{C}$,
\begin{equation}\label{2.4}
\Psi (y,\lambda p) = \lambda \Psi (y,p).
\end{equation}
\end{lemm}
{\bf Preuve} Montrons d'abord (\ref{2.3}). 
Fixons $\lambda \in \Bbb{C}$ et d\'erivons la relation (\ref{1.6}) par
rapport \`a $\overline{z}^j$.
Il vient:

$$\overline{\lambda}{\partial F\over \partial \overline{z}^j}(y,\lambda z) =
{\partial F(y,\lambda z)\over \partial \overline{z}^j} =
|\lambda |^2{\partial F\over \partial \overline{z}^j}(y,z).$$
En simplifiant par $\overline{\lambda}$, cela donne (\ref{2.3}). A pr\'esent,
puisque $\Psi$ est l'inverse de $2{\partial F\over \partial \overline{z}}$,
on a, notant $z = \Psi(y,p)$,

$$\Psi(y,\lambda p) = \Psi(y,\lambda 2{\partial F\over \partial
\overline{z}}(y,z)) =
\Psi(y,2{\partial F\over \partial \overline{z}}(y,\lambda z)) = \lambda z =
\lambda \Psi(y,p),$$
ce qui donne (\ref{2.4}). {\em CQFD}.

\section{G\'eom\'etrie $\Bbb{C}$-finslerienne}
La g\'eom\'etrie finslerienne est obtenue en d\'efinissant une fonction $F$
sur le fibr\'e tangent $T{\cal N}$, \`a valeurs
dans les r\'eels positifs, qui est homog\`ene de degr\'e 1, c'est \`a dire
$F(y,\lambda v)=|\lambda |F(y,v)$, pour tout
$(y,v)\in T{\cal N}$ et tout $\lambda \in \Bbb{R}$. Elle consiste en une
vision g\'eom\'etrique des probl\`emes
variationnels \`a une variable invariants par les diff\'eomorphismes de
$\Bbb{R}$ (cf [Ch]).
La g\'eom\'etrie $\Bbb{C}$-finslerienne, que nous allons d\'ecrire plus bas, mod\'elise de
fa\c{c}on g\'eom\'etrique les probl\`emes variationnels \`a deux variables invariants par 
transformation conforme et se pr\'esente
comme une version complexe de le g\'eom\'etrie finslerienne. Par exemple, une belle construction de
la g\'eom\'etrie finslerienne consiste \`a d\'efinir sur le fibr\'e tangent
$T{\cal N}$ un tenseur m\'etrique
(un produit scalaire sur chaque $T_y{\cal N}$) $g_{ij}(y,v)={1\over
2}(\partial ^2/\partial v^i\partial v^j)F^2(y,v)$. Ce produit
scalaire est homog\`ene de degr\'e z\'ero, c'est \`a dire $g_{ij}(y,\lambda
v) = \lambda g_{ij}(y,v)$. Cela signifie qu'il est d\'efini
sur le fibr\'e projectif tangent $PT{\cal N}$. On fabrique ainsi une
m\'etrique sur chaque fibre du fibr\'e
$P{\cal F}_1$, dont la vari\'et\'e base est $PT{\cal N}$ et la fibre en un
point $(y,\Bbb{R}v)$ est $T_y{\cal N}$.
Dans ce qui suit, nous pr\'esentons une construction similaire en
g\'eom\'etrie $\Bbb{C}$-finslerienne.

\subsection{Tenseurs m\'etriques}
Nous diff\'erentions la relation (\ref{2.3}) une fois de plus, par rapport aux
variables $z$.
Introduisons la notation

$${\partial \over \partial z^k} = {1\over 2}\left(
{\partial \over \partial z_1^k} - i{\partial \over \partial z_2^k} \right) ,$$
et appliquons cet op\'erateur \`a (\ref{2.3}), Il vient

$${\lambda}{\partial ^2F\over \partial \overline{z}^j\partial
z^k}(y,\lambda z) =
\lambda {\partial ^2F\over \partial \overline{z}^j\partial z^k}(y,z).$$
En simplifiant par $\lambda$ et en identifiant les parties r\'eelles et
imaginaires, nous trouvons

\begin{equation}\label{4.1}\begin{array}{c}
G_{jk}^{11}(y,\lambda z) + G_{jk}^{22}(y,\lambda z) = G_{jk}^{11}(y,z) +
G_{jk}^{22}(y,z)\\
\\
G_{jk}^{12}(y,\lambda z) - G_{jk}^{21}(y,\lambda z) = G_{jk}^{12}(y,z) -
G_{jk}^{21}(y,z).
\end{array}
\end{equation}
Nous sommes ainsi conduits \`a poser
\begin{equation}\label{4.2}\begin{array}{c}
g_{jk}(y,z) = {1\over 2}\left( G_{jk}^{11}(y,z) + G_{jk}^{22}(y,z)\right) \\
\\
\omega_{jk}(y,z) = {1\over 2}\left( G_{jk}^{12}(y,z) -
G_{jk}^{21}(y,z)\right) .
\end{array}
\end{equation}
Il est imm\'ediat que $g_{jk}$ est sym\'etrique, $\omega_{jk}$ est
antisym\'etrique et que l'on a
$g_{jk}(y,z)v^jv^k\geq c|v|^2$. A partir de $g_{jk}$ et de $\omega_{jk}$,
on peut former le
tenseur $h_{jk}:=g_{jk} - i\omega_{jk}$, qui satisfait $\overline{h_{kj}} =
h_{jk}$,
c'est \`a dire, qui est hermitien.
Notons ${\cal F}$ le fibr\'e image inverse de $T{\cal N}$ par
la fibration $T{\cal N}^{\Bbb{C}}\longrightarrow {\cal N}$.
Nous avons obtenu le r\'esultat suivant. 

\begin{lemm}
{\bf et d\'efinition} Les tenseurs $g$ et $\omega$ donn\'es par (\ref{4.2})
d\'efinissent une
``m\'etrique hermitienne''

$$h_{jk} = g_{jk}-i\omega_{jk} = 2{\partial ^2F\over \partial \overline{z^j}\partial z^k}\hbox{ sur }{\cal F}.$$
Cette m\'etrique est homog\`ene complexe de degr\'e z\'ero sur
chaque $T_y{\cal N}^{\Bbb{C}}$.    
En d'autres termes, nous introduisons $PT{\cal N}^{\Bbb{C}}$, le fibr\'e
dont la fibre en
$y$ est l'espace projectif complexe $PT_y{\cal N}^{\Bbb{C}}$ (le quotient
de $T{\cal N}^{\Bbb{C}}$ par $\Bbb{C}$).
Nous notons $P{\cal F}$ le fibr\'e image inverse de $T{\cal N}$ par
la projection $PT{\cal N}^{\Bbb{C}}\longrightarrow {\cal N}$. Alors le
tenseur $h$ d\'efinit une m\'etrique hermitienne
sur $P{\cal F}$.
\end{lemm}
Nous allons \`a pr\'esent voir qu'il existe des relations entre les
``vitesses g\'en\'eralis\'ees'' $z\in T_y{\cal N}^{\C}$, les
impulsions g\'en\'eralis\'ees $p\in T_y^{\star}{\cal N}^{\C}$ et la
m\'etrique hermitienne $g-i\omega$ similaires \`a celles de la
g\'eom\'etrie riemannienne. Repartons de (\ref{1.5b}), que nous pouvons
r\'e\'ecrire sous la forme

\begin{equation}\label{4.3}
{1\over 2}\overline{p}_jz^j = {\partial F\over \partial z^j}(y,z)z^j = F(y,z).
\end{equation}

D\'erivons (\ref{4.3}) par rapport \`a $\overline{z}^k$:

$${\partial ^2F\over \partial z^j\partial \overline{z}^k}(y,z)z^j =
{\partial F\over \partial \overline{z}^k}(y,z).$$
Remarquant que ${\partial ^2F\over \partial z^j\partial \overline{z}^k} =
{1\over 2}(g_{jk}+i\omega_{jk})$, nous en d\'eduisons

\begin{equation}\label{4.5}
(g_{jk}+i\omega_{jk})(z^j_1+iz^j_2) = 2{\partial F\over \partial
\overline{z}^k}(y,z) = (p^1_k +ip^2_k).
\end{equation}
Enfin, utilisant (\ref{4.3}) et (\ref{4.5}), on obtient:

\begin{equation}\label{4.6}
\begin{array}{ccl}
F(y,z) & = & {1\over 2}\left( p^1_k -ip^2_k\right) \left(
z^k_1+iz^k_2\right) \\
& = & {1\over 2}\left( g_{jk}-i\omega_{jk}\right) (y,z)\left(
z^j_1-iz^j_2\right) \left( z^k_1+iz^k_2\right) \\
& = & {1\over 2}\left( g_{jk}-i\omega_{jk}\right) (y,z)\left( (z^j_1z^k_1 +
z^j_2z^k_2) +i(z^j_1z^k_2 -z^j_2z^k_1)\right)\\

& = & {1\over 2}\left( g_{jk}(y,z)\left( z^j_1z^k_1 + z^j_2z^k_2\right)  +
\omega_{jk}(y,z)\left( z^j_1z^k_2 -z^j_2z^k_1\right) \right) .
\end{array}
\end{equation}

\begin{rema}
On a aussi $F(y,z)= {1\over 2}h_{jk}(y,z)\overline{z^j}z^k$.
\end{rema}
Il est maintenant possible de reformuler l'action d'une application
$u:\Sigma \longrightarrow {\cal N}$ de la fa\c{c}on
suivante. Nous supposons que $\overline{\partial}u(t):={1\over 2}\left(
{\partial u\over \partial t^1} + i {\partial u\over \partial t^2}\right)$
ne s'annule nulle part. Nous d\'efinissons l'unique rel\`evement $Pu:\Sigma
\longrightarrow PT{\cal N}^{\C}$ de $u$ qui soit
{\em horizontal}, c'est \`a dire:

\begin{equation}\label{4.7}
\forall t\in \Sigma ,\ Pu(t) = (u(t), \Bbb{C}\overline{\partial}u(t)),
\end{equation}
o\`u $\Bbb{C}\overline{\partial}u(t)$ est la droite complexe dans
$T_{u(t)}{\cal N}^{\C}$ engendr\'ee par $\overline{\partial}u(t)$.
Alors

\begin{lemm}
L'action de $u:\Sigma \longrightarrow {\cal N}$ est \'egale \`a l'int\'egrale

$${\cal E}[u] = \int _{\Sigma} {1\over 2} g_{jk}(Pu)(\partial _1u^j\partial
_1u^k + \partial _2u^j\partial _2u^k)dt^1\wedge dt^2
+ {1\over 2} \omega_{jk}(Pu)du^j\wedge du^k.$$
\end{lemm}
Dans cette derni\`ere formule, l'action appara\^\i t comme la somme d'une
int\'egrale de Dirichlet et d'une int\'egrale
$\int _{\Sigma} Pu^{\star}\omega$, o\`u $\omega = {1\over
2}\omega_{jk}(y,\Bbb{C}z)dy^j\wedge dy^k$. Bien \'evidemment, Il faut faire
attention que tout cela n'est valable que lorsque la condition
d'horizontalit\'e (\ref {4.7}) sur $Pu$ est v\'erifi\'ee.

\subsection{Equation d'Euler-Lagrange}
Ecrivons les \'equations v\'erifi\'ees par les points critiques dans notre
formalisme. Nous partons de

$${\partial \over \partial t^{\alpha}}\left( {\partial F\over \partial
z^j_{\alpha}}(u,2\overline{\partial}u)\right)
={\partial F\over \partial y^j}(u,2\overline{\partial}u),$$
qui donne en d\'eveloppant:

\begin{equation}\label{4.8}
{\partial ^2F\over \partial z^j_{\alpha}\partial z^k_{\beta}}{\partial
^2u^k\over \partial t^{\alpha}\partial t^{\beta}} +
{\partial ^2F\over \partial z^j_{\alpha}\partial y^k}{\partial u^k\over
\partial t^{\alpha}} -
{\partial F\over \partial y^j} = 0.
\end{equation}
Pour interpr\'eter cette relation, nous d\'erivons les relations
(\ref{4.5}) et (\ref{4.6}) par rapport à $y$.
Pour (\ref{4.6}) ou $F(y,z)={1\over 2}h_{kl}(y,z)\overline{z^k}z^l$, on
obtient:

\begin{equation}\label{4.9}
{\partial F\over \partial y^j} = {1\over 2} {\partial h_{kl}\over \partial
y^j}\overline{z^k}z^l.
\end{equation}
Et pour (\ref{4.5}) ou $2{\partial F\over \partial \overline{z^j}} =
h_{jl}z^l$,

$$2{\partial ^2F\over \partial \overline{z^j}\partial y^k} =
{\partial h_{jl}\over \partial y^k}z^l,$$
relation qui entraîne, notant $z^k={\partial u^k\over \partial t^1}
+i{\partial u^k\over \partial t^2}$,

\begin{equation}\label{4.10}
{\partial ^2F\over \partial z^j_{\alpha}\partial y^k}{\partial u^k\over
\partial t^{\alpha}} =
2 Re \left( {\partial ^2F\over \partial \overline{z^j}\partial
y^k}\overline{z^k}\right) =
Re\left( {\partial h_{jl}\over \partial y^k}z^l\overline{z^k}\right) =
{1\over 2} {\partial h_{jl}\over \partial y^k}z^l\overline{z^k} + {1\over
2} {\partial h_{lj}\over \partial y^k}\overline{z^l}z^k.
\end{equation}
Donc, utilisant (\ref{4.9}) et (\ref{4.10}), on obtient

\begin{equation}\begin{array}{ccl}
\displaystyle {\partial ^2F\over \partial z^j_{\alpha}\partial
y^k}{\partial u^k\over \partial t^{\alpha}} - {\partial F\over \partial
y^j} & = &
\displaystyle {1\over 2} \left[ {\partial h_{jl}\over \partial
y^k}z^l\overline{z^k} + {\partial h_{lj}\over \partial y^k}\overline{z^l}z^k
- {\partial h_{kl}\over \partial y^j}\overline{z^k}z^l\right] \\
& = & \displaystyle {1\over 2} \left[ {\partial g_{jl}\over \partial y^k} +
{\partial g_{kj}\over \partial y^l}
- {\partial g_{kl}\over \partial y^j}\right] (z_1^kz_1^l+z_2^kz_2^l) \\
&  & + \displaystyle {1\over 2} \left[ {\partial \omega_{jl}\over \partial
y^k} + {\partial \omega_{kj}\over \partial y^l}
+ {\partial \omega_{lk}\over \partial y^j}\right] (z_1^kz_2^l-z_2^kz_1^l).
\end{array}
\end{equation}
Nous en d\'eduisons que l'\'equation (\ref{4.8}) peut s'\'ecrire

\begin{equation}\label{4.11}
g^{jm}G^{\alpha \beta}_{jk}{\partial ^2u^k\over \partial t^{\alpha}\partial
t^{\beta}}
+ \Gamma ^m_{kl}(z_1^kz_1^l+z_2^kz_2^l) - {1\over
2}g^{jm}(d\omega)_{jkl}(z_1^kz_2^l-z_2^kz_1^l) = 0,
\end{equation}
où 

$$
\Gamma ^m_{kl} = {1\over 2} g^{jm}\left( {\partial g_{jl}\over \partial
y^k} + {\partial g_{kj}\over \partial y^l}
- {\partial g_{kl}\over \partial y^j}\right) 
$$
et
$$
(d\omega)_{jlk} = {\partial \omega_{jl}\over \partial y^k} + {\partial
\omega_{kj}\over \partial y^l}
+ {\partial \omega_{lk}\over \partial y^j}.
$$

Enfin, il est possible d'expliciter diff\'eremment le terme d'ordre 2 dans
l'\'equation d'Euler-Lagrange, en introduisant
les tenseurs sym\'etriques $a_{jk}$ et $b_{jk}$ donn\'es par

\begin{equation}\label{4.12}
a_{jk}-ib_{jk} := {1\over 2}\left( G^{11}_{jk}-G^{22}_{jk}\right) - {i\over
2}\left( G^{12}_{jk}+G^{21}_{jk}\right) 
= 2{\partial ^2F\over \partial z^j\partial z^k},
\end{equation}
tels que

$$
G^{11}_{jk} = g_{jk} + a_{jk},\ G^{12}_{jk} = \omega_{jk} + b_{jk},\
G^{21}_{jk} = - \omega_{jk} + b_{jk} ,\ G^{22}_{jk} = g_{jk} - a_{jk}.
$$
En effet, on a alors

$$g^{jm}G^{\alpha \beta}_{jk}{\partial ^2u^k\over \partial
t^{\alpha}\partial t^{\beta}} =
\Delta u^m + g^{jm}a_{jk}(\partial _1^2 -\partial _2^2)u^k +
2g^{jm}b_{jk}\partial _1\partial _2u^k,$$
où $\partial _1 ={\partial \over \partial t^1}$, $\partial _2 ={\partial
\over \partial t^2}$.
Donc l'\'equation d'Euler est dans ces notations

\begin{equation}\label{4.13}
\Delta u^m + g^{jm}a_{jk}(\partial _1^2 -\partial _2^2)u^k +
2g^{jm}b_{jk}\partial _1\partial _2u^k
+ \Gamma ^m_{kl}(z_1^kz_1^l+z_2^kz_2^l) - {1\over
2}g^{jm}(d\omega)_{jkl}(z_1^kz_2^l-z_2^kz_1^l) = 0.
\end{equation}

Cette \'equation serait identique à l'\'equation d'une application
harmonique à valeurs dans ${\cal N}$ munie de la m\'etrique
$g_{jk}$, en pr\'esence d'une forme $\omega_{jk}$, si les termes en
$a_{jk}$ et $b_{jk}$ \'etaient nuls.
\begin{lemm}
Les tenseurs $a_{jk}$ et $b_{jk}$ sont nuls si et seulement si $F(y, z) =
{1\over 2}h_{jk}(y)\overline{z^j}z^k$.
\end{lemm}
{\bf Preuve} D'apr\`es (\ref{4.12}), $a_{jk} = b_{jk} = 0$ si et seulement si
${\partial ^2F\over \partial z^j\partial z^k} = 0$.
Comme $F$ est r\'eel, cette relation entra\^{\i}ne aussi ${\partial ^2F\over
\partial \overline{z^j}\partial \overline{z^k}} = 0$
et donc n\'ecessairement $F$ est de la forme $F(y,z,\overline{z}) = \sum
_{j=1}^n A_j(y,z)\overline{z^j} + A_0(y,z)
= \sum _{j=1}^n B_j(y,\overline{z})z^j + B_0(y,\overline{z})$.
De plus,  $F$ doit \^etre
homog\`ene de degr\'e 2 en $z^j_{\alpha}$, et donc $F(y, z) = {1\over
2}\eta_{jk}(y)\overline{z^j}z^k$, o\`u $\eta_{jk}$ est un tenseur
hermitien. Les relations (\ref{4.2}) permettent alors de conclure que
$\eta_{jk} = h_{jk}$
\footnote{Remarquons que, bien que $a_{jk}$ et $b_{jk}$ soient non nuls en
g\'en\'eral, on a toujours $\left( a_{jk} -ib_{jk}\right) z^j = 0$.
En effet, en d\'erivant (\ref{4.3}) par rapport \`a $z^k$, on obtient

$${\partial ^2F\over \partial z^j\partial z^k}(y,z)z^j + {\partial F\over
\partial z^k}(y,z)= {\partial F\over \partial z^k}(y,z),$$
d'o\`u, en simplifiant ${\partial ^2F\over \partial z^j\partial
z^k}(y,z)z^j = 0$, qui donne cette relation d'après (\ref{4.12})}. {\em CQFD}.

\section{Approches hamiltoniennes}
Dans le formalisme hamiltonien classique pour des probl\`emes variationnels
\`a une variable, lorsque l'hypoth\`ese de
Legendre est satisfaite, on remplace les
variables $t\in \Bbb{R}$ et $(y,z)\in T{\cal N}$ par $t\in \Bbb{R}$ et
$(y,p)\in T^{\star}{\cal N}$. L'\'equivalence
entre les deux syst\`emes de coordonn\'ees repose sur le fait que
$p_i={\partial L\over \partial z^i}(y,z)$ est un diff\'eomorphisme
entre $T{\cal N}$ et $T^{\star}{\cal N}$. On d\'efinit l'hamiltonien sur
$\Bbb{R}\times T^{\star}{\cal N}$
par $H(t,y,{\partial L\over \partial z}(t,y,z)) = {\partial L\over \partial
z^i}(t,y,z)z^i - L(t,z,y)$.
Alors $H$ joue deux r\^oles: d'une part, lorsque le probl\`eme est
ind\'ependant du temps, il est une quantit\'e conserv\'ee,
d'autre part, il dicte la dynamique du probl\`eme par les \'equations de
Hamilton. On peut m\^eme adopter une description plus
sym\'etrique en espace et en temps en consid\'erant sur $T^{\star}\left(
\Bbb{R}\times {\cal N}\right)$ l'hamiltonien
$\hat{H}(t,y^1,...,y^n,p_0,p_1,...,p_n) = H(t,y,p) +p_0$. Une fa\c{c}on
d'obtenir les \'equations de Hamilton est d'\'ecrire
que la trajectoire $\gamma (t) =
(t,y^1(t),...,y^n(t),p_0(t),p_1(t),...,p_n(t))$ est un point critique de la
fonctionnelle

$${\cal A} = \int _I p_i(t)dy^i(t) - H(t,y(t),p(t))dt$$
ou encore de
$${\cal A} = \int _I p_idy^i + p_0dt,$$
avec la contrainte $\hat{H}(t,y,p_0,p) = 0$. Dans le cas de cette derni\`ere formulation,
c'est la contrainte qui dicte la dynamique.

L'id\'eal serait de pouvoir faire de m\^eme pour le calcul des variations
\`a plusieurs variables. Malheureusement les choses ne marchent pas aussi bien. Des
constructions partiellement satisfaisantes
existent. Nous en pr\'esentons deux ici, la th\'eorie de Weyl et celle de
Carath\'eodory. Mais il y en a en r\'ealit\'e
une infinit\'e (th\'eories de de Donder, de Boerner...). De plus, les deux
r\^oles jou\'es par l'hamiltonien pour les probl\`emes
\`a une variable sont en g\'en\'eral r\'epartis entre deux objets
diff\'erents: un hamiltonien scalaire, dont le r\^ole est dynamique et
un hamiltonien tensoriel (le tenseur \'energie-impulsion), qui contient les
quantit\'es conserv\'ees. Dans la suite, nous exposerons
les th\'eories de Weyl et de Carath\'eodory de fa\c{c}on succincte et
axiomatique.
Pour une pr\'esentation plus approfondie et les liens avec la th\'eorie
d'Hamilton-Jacobi, nous recommandons de lire l'ouvrage de H. Rund [Ru]
\footnote{voir \'egalement dans ce livre les r\'ef\'erences bibliographiques des auteurs suivants: E. T. Davies, A. Kawaguchi, A. Kawaguchi et Y. Katsurada, A. Kawaguchi et K. Tandai, L. Berwald}.

\subsection{Th\'eorie de Weyl}
Nous avons juste besoin ici de la condition de Legendre globale. Comme nous
l'avons vu au paragraphe 3.2, cela garantit l'existence d'une application
inverse de
$z\mapsto {\partial F\over \partial z^j_{\alpha}}(y,z)$, que nous avons
not\'ee $\Psi = \psi_1+i\psi_2$.
Nous d\'efinissons alors l'hamiltonien de Weyl $H:T^{\star}{\cal
N}^{\Bbb{C}}\longrightarrow \Bbb{R}$ par

\begin{equation}\label{weyl1}
H(y,p) := p^{\alpha}_j\psi^j_{\alpha}(y,p) - F\left( y,\Psi(y,p)\right) .
\end{equation}
Un calcul direct donne
\begin{equation}\label{weyl2}
{\partial H\over \partial p^{\alpha}_j}(y,p) = \psi^j_{\alpha}(y,p)
\end{equation}
et
\begin{equation}\label{weyl3}
{\partial H\over \partial u^j}(y,p) = -{\partial F\over \partial u^j}\left(
y,\Psi(y,p)\right) .
\end{equation}

Nous remarquons qu'en substituant $p^{\alpha}_j={\partial F\over \partial
z^j_{\alpha}}(y,z)$,
avec $z=\Psi(p)$ dans (\ref{weyl1}), nous obtenons

$$H\left( y,{\partial F\over \partial z}(y,z)\right) =
{\partial F\over \partial z^j_{\alpha}}(y,z)z^j_{\alpha} -F(y,z) = F(y,z),$$
en vertu de (\ref{1.5b}). Et donc
\begin{equation}\label{weyl4}
H(y,p) = F(y, \Psi(y,p)).
\end{equation}
Gr\^ace \`a cette identit\'e, nous d\'eduisons imm\'ediatement de
l'hypoth\`ese (\ref{1.6})
et du Lemme 1 que $\forall (y,p)\in T^{\star}{\cal N}^{\Bbb{C}}$, $\forall
\lambda \in \Bbb{C}$,
$H(y, \lambda p) = |\lambda |^2H(y,p)$.

Egalement, en utilisant les relations (\ref{4.5}) et (\ref{4.6}), on
obtient que
$H(y,p) = {1\over 2}\eta ^{jk}(y,p)\overline{p}_jp_k$, o\`u $\eta
^{jk}(y,p)$ est
d\'efini par $\eta ^{jk}(y,p)h_{kl}(y,z) = \delta ^j_l$.

De plus, le syst\`eme d'\'equations d'Euler-Lagrange
${\partial \over\partial t^{\alpha}}\left( {\partial L\over \partial
z^j_{\alpha}}(u,du)\right)
={\partial L\over \partial y^j}(u,du)$ peut \^etre transform\'e en y
ajoutant les relations
(\ref{weyl2}), reliant ${\partial u^j\over \partial t^{\alpha}}$ \`a
$p^{\alpha}_j$.
En utilisant aussi (\ref{weyl3}), cela donne les \'equations de Hamilton
g\'en\'eralis\'ees, comme suit

\begin{equation}\label{weyl5}
\left\{
\begin{array}{ccc}
\displaystyle {\partial u^j\over \partial t^{\alpha}} & = & \displaystyle
{\partial H\over \partial p_j^{\alpha}}(u,p)\\
\displaystyle {\partial p^1_j\over \partial t^1} + {\partial p^2_j\over
\partial t^2}& = & \displaystyle -{\partial H\over \partial y^j}(u,p).
\end{array}\right.
\end{equation}
Nous remarquons que les premi\`eres relations entra\^{\i}nent aussi les
relations de compatibilit\'e

$${\partial \over \partial t^1}\left( {\partial H\over \partial
p^2_j}(u,p)\right) -
{\partial \over \partial t^2}\left( {\partial H\over \partial
p^1_j}(u,p)\right) =0.$$

De plus, pour un ouvert $\Omega\subset \Bbb{C}$, $(u,p):\Omega
\longrightarrow T{\cal N}^{\Bbb{C}}$ est  solution de (\ref{weyl5})
si et seulement si son graphe $\Gamma =\{(t,u(t),p(t))\subset \Omega \times
T^{\star}{\cal N}^{\Bbb{C}}/t\in \Omega \}$ est un point critique de

$$\int _{\Gamma} p^1_jdu^j\wedge dt^2 + p^2_jdt^1\wedge du^j -
H(t,u,p)dt^1\wedge dt^2.$$

\subsection{Th\'eorie de Carath\'eodory-Rund}
Dans ce qui suit, nous reprenons et d\'eveloppons 
le formalisme canonique construit par H. Rund [Ru], [Ru1], dans le
but d'interpr\'eter les \'equations d'Hamilton-Jacobi de Carath\'eodory
(voir le paragraphe 5.3).
L'id\'ee est d'utiliser comme variables d'impulsion des quantit\'es
$\epsilon_{\beta}^{\alpha}$ et $\pi_j^{\alpha}$, pour $\alpha , \beta =1,2$
et $j=1,...,n$, telles que

\begin{equation}\label{cara1}
F(y,z) + w = \left| \begin{array}{cc}
\pi^1_jz^j_1+\epsilon_1^1 & \pi^1_kz^k_2+\epsilon_2^1\\
\pi^2_jz^j_1+\epsilon_1^2 & \pi^2_kz^k_2+\epsilon_2^2
\end{array}\right| ,
\end{equation}
et qu'il existe une matrice $2\times 2$ inversible $T$, telle que

\begin{equation}\label{cara2}
\pi = T {\partial F\over \partial z}\hbox{ et }
\epsilon = wT - T H,
\end{equation}
avec

$$\pi = \left( \begin{array}{ccc}
\pi_1^1 & \dots & \pi_n^1\\
\pi_1^2 & \dots & \pi_n^2\\
\end{array}\right) ,
\epsilon = \left( \begin{array}{cc}
\epsilon_1^1 & \epsilon_2^1\\
\epsilon_1^2 & \epsilon_2^2\\
\end{array}\right) ,
$$
et
$${\partial F\over \partial z} = \left( \begin{array}{ccc}
{\partial F\over \partial z^1_1} & \dots & {\partial F\over \partial z^n_1}\\
{\partial F\over \partial z^1_2} & \dots & {\partial F\over \partial z^n_2}
\end{array}\right) ,
H = \left( \begin{array}{cc}
H^1_1 & H^2_1\\ H^1_2 & H^2_2
\end{array}\right) .
$$
Ici, $w$ est un param\`etre r\'eel (qui n'appara\^{\i}t pas dans [Ru]),
dont le r\^ole
sera \'eclairci plus tard.
Remarquons que si l'on note

$$Z = \left( \begin{array}{cc}
z^1_1 & z^1_2\\
\vdots & \vdots \\
z^n_1 & z^n_2
\end{array}\right) ,$$
on a $H = {\partial F\over \partial z}Z-F\1 _2$ et donc (\ref{cara1})
entra\^{\i}ne que

$$F(y,z) + w = \det (\pi Z + \epsilon ) =
\det \left( T {\partial F\over \partial z}Z + wT -T H\right) =
\det ((F+w)T )= (F+w)^2 \det T.$$

Nous en d\'eduisons que si $F(y,z) + w \neq 0$,

\begin{equation}\label{cara3}
(F(y,z) + w) \det T = 1.
\end{equation}
R\'eciproquement, il est imm\'ediat que pour tout $T$ qui satisfait
(\ref{cara3}),
si $\epsilon$ et $\pi$ sont d\'efinis par (\ref{cara2}), alors
(\ref{cara1}) est v\'erifi\'e.
En excluant le cas $F(y,z)+w=0$, on en d\'eduit qu'\'etant donn\'e $w$, \`a
tout $z$,
on peut associer un couple $(\epsilon,\pi)$ satisfaisant
(\ref{cara1}) et (\ref{cara2}). Cette solution $(\epsilon,\pi)$ n'est pas unique, puisque l'on peut changer $T$ en $gT$
(ce qui revient \`a changer $(\epsilon,\pi)$ en $(g\epsilon,g\pi)$),
pour tout $g\in SL(2, \Bbb{R})$. Donc, dans les variables $(\epsilon,\pi)$,
l'ensemble des quantit\'es ``observables'' \footnote{c'est \`a dire les fonctions de $(y,z)$} est invariant par le groupe de jauge $SL(2, \Bbb{R})$.

\subsubsection{Correspondance de Legendre-Carath\'eodory}
Nous allons v\'erifier que l'on peut remplacer les variables
$(y,z,w)$ par $(y,\epsilon,\pi)$. Dans ce qui suit, nous supposerons que
$F(y,z) + w\neq 0$.
Ce qui pr\'ec\`ede montre qu'\'etant donn\'e $z,w$, l'on peut toujours trouver
$(\epsilon,\pi)$ tels que (\ref{cara1}) et (\ref{cara2}) aient lieu. La
r\'eciproque
est un peu plus d\'elicate. Commen\c{c}ons par caract\'eriser $T$ en fonction
de $z$, $\epsilon$ et $\pi$.

\begin{lemm}
Si $T$ est solution de (\ref{cara1}) et (\ref{cara2}), alors 

$$T = {\pi Z + \epsilon \over \det (\pi Z + \epsilon )},$$
ou, de fa\c{c}on \'equivalente, $T ^{-1}$ est la comatrice  de $\pi Z + \epsilon$,
c'est \`a dire,

\begin{equation}\label{cara4}
T ^{-1} = \left( \begin{array}{cc}
\pi^2_kz^k_2+\epsilon_2^2 & -(\pi^1_kz^k_2+\epsilon_2^1)\\
-(\pi^2_jz^j_1+\epsilon_1^2) & \pi^1_jz^j_1+\epsilon_1^1
\end{array}\right)
:= P = \left( \begin{array}{cc}
P^1_1 & P^1_2 \\ P^2_1 & P^2_2\end{array}\right) .
\end{equation}
\end{lemm}
{\bf Preuve} A partir de la relation $F\1 _2 = {\partial F\over \partial z}Z-H$, nous 
obtenons

$$(F(y,z) + w)T =  T {\partial F\over \partial z}Z
- T H + wT = \pi Z + \epsilon ,$$
d'o\`u, utilisant (\ref{cara3}), ${T \over det T } = \pi Z + \epsilon$.
Cela implique $(\det T )^{-1} = \det (\pi Z + \epsilon)$, donc
$T = {\pi Z + \epsilon \over \det (\pi Z + \epsilon )}$. La relation
(\ref{cara4}) s'ensuit. {\em CQFD}.\\

A pr\'esent, pour trouver $z$ en fonction de $\epsilon$ et
$\pi$, il
suffit de r\'esoudre (\ref{cara2}) en y substituant $T$ selon (\ref{cara4}),
c'est \`a dire r\'esoudre le syst\`eme

\begin{equation}\label{cara5}
{\partial F\over \partial z^j_1}(y,z) = \left| \begin{array}{cc}
\pi^1_j & \pi^1_kz^k_2+\epsilon_2^1\\
\pi^2_j & \pi^2_kz^k_2+\epsilon_2^2
\end{array}\right| ,\ 
{\partial F\over \partial z^k_2}(y,z) = \left| \begin{array}{cc}
\pi^1_jz^j_1+\epsilon_1^1 & \pi^1_k\\
\pi^2_jz^j_1+\epsilon_1^2 & \pi^2_k
\end{array}\right| ,
\end{equation}

\begin{equation}\label{cara6}
\begin{array}{c}\displaystyle
H^1_2(y,z) = - \left| \begin{array}{cc}
\epsilon_2^1 & \pi^1_kz^k_2+\epsilon_2^1\\
\epsilon_2^2 & \pi^2_kz^k_2+\epsilon_2^2
\end{array}\right| 
= \left| \begin{array}{cc}
\pi^1_jz^j_2 & \pi^1_kz^k_2+\epsilon_2^1\\
\pi^2_jz^j_2 & \pi^2_kz^k_2+\epsilon_2^2
\end{array}\right| ,\\
\displaystyle
H^2_1(y,z) = - \left| \begin{array}{cc}
\pi^1_jz^j_1+\epsilon_1^1 & \epsilon_1^1\\
\pi^2_jz^j_1+\epsilon_1^2 & \epsilon_1^2
\end{array}\right| 
= \left| \begin{array}{cc}
\pi^1_jz^j_1+\epsilon_1^1 & \pi^1_kz^k_1\\
\pi^2_jz^j_1+\epsilon_1^2 & \pi^2_kz^k_1
\end{array}\right| .
\end{array}
\end{equation}
\noindent et
\begin{equation}\label{cara7}
H^1_1(y,z) = w - \left| \begin{array}{cc}
\epsilon_1^1 & \pi^1_kz^k_2+\epsilon_2^1\\
\epsilon_1^2 & \pi^2_kz^k_2+\epsilon_2^2
\end{array}\right| ,\ 
H^2_2(y,z) = w -\left| \begin{array}{cc}
\pi^1_jz^j_1+\epsilon_1^1 & \epsilon_2^1\\
\pi^2_jz^j_1+\epsilon_1^2 & \epsilon_2^2
\end{array}\right| .
\end{equation}

Nous allons voir que,
\begin{itemize}
\item sous des hypoth\`eses g\'en\'eriques sur $\epsilon$ et $\pi$,
le syst\`eme (\ref{cara5}) admet une unique solution $z$
\item on peut v\'erifier qu'alors $z$ est automatiquement solution
de (\ref{cara6})
\item enfin, il existe un unique $w$ tel que $z$ et $w$ soient solutions de
(\ref{cara7}).
\end{itemize}

Auparavant, nous remarquons sur (\ref{cara5}), (\ref{cara6}) et (\ref{cara7})
que $z$ et $w$ ne d\'ependent que des coefficients

$$A_{j,k}:= \left| \begin{array}{cc}
\pi^1_j & \pi^1_k\\
\pi^2_j & \pi^2_k
\end{array}\right| ,\
A_{j,n+2}:=\left| \begin{array}{cc}
\pi^1_j & \epsilon_2^1\\
\pi^2_j & \epsilon_2^2
\end{array}\right| ,
A_{n+1,k}:=\left| \begin{array}{cc}
\epsilon_1^1 & \pi^1_k\\
\epsilon_1^2 & \pi^2_k
\end{array}\right| \hbox{ et }
A_{n+1,n+2}:=\left| \begin{array}{cc}
\epsilon_1^1& \epsilon_2^1\\
\epsilon_1^1& \epsilon_2^2
\end{array}\right| ,$$
quantit\'es invariantes par un changement de $(\epsilon,\pi)$ en $(g\epsilon,g\pi)$,
pour $g\in SL(2, \Bbb{R})$. C'est pourquoi nous adopterons plut\^ot
ces notations ($\left( A_{jk}\right) _{1\leq j,k\leq n+2}$ formant une
matrice antisym\'etrique
$(n+2)\times (n+2)$) et les \'equations (\ref{cara5}), (\ref{cara6}),
(\ref{cara7}) se r\'e\'ecrivent respectivement

\begin{equation}\label{cara8}
{\partial F\over \partial z^j_1}(y,z) = A_{j,k}z_2^k + A_{j,n+2} ,\ 
{\partial F\over \partial z^k_2}(y,z) = A_{j,k}z_1^j + A_{n+1,k} ,
\end{equation}

\begin{equation}\label{cara9}
H^1_2(y,z) =  A_{j,n+2}z_2^j ,\ 
H^2_1(y,z) =  A_{n+1,k}z_1^k 
\end{equation}
et

\begin{equation}\label{cara10}
H^1_1(y,z) =  w - \left( A_{n+1,k}z_2^k + A_{n+1,n+2}\right) ,\ 
H^2_2(y,z) =  w - \left( A_{j,n+2}z_1^j + A_{n+1,n+2}\right) .
\end{equation}

\noindent {\bf R\'esolution de (\ref{cara8}).} Nous consid\'erons la
fonctionnelle

$$
W(y,z,\epsilon,\pi) := \left| \begin{array}{cc}
\pi^1_jz^j_1+\epsilon_1^1 & \pi^1_kz^k_2+\epsilon_2^1\\
\pi^2_jz^j_1+\epsilon_1^2 & \pi^2_kz^k_2+\epsilon_2^2
\end{array}\right| - F(y,z)
$$
ou
$$ W(y,z,A):= 
\left( \sum _{j,k=1}^n A_{j,k}z_1^jz_2^k + \sum _{j=1}^n A_{j,n+2}z_1^j
+\sum _{k=1}^n A_{n+1,k}z_2^k + A_{n+1,n+2} \right) - F(y,z) .$$
Toute solution $z$ de (\ref{cara8}) est obtenue en fixant
$A:=\left( A_{jk}\right) _{1\leq j,k\leq n+2}$ et
en extr\'emisant $z\longmapsto W(y,z,A)$. En effet, il est
imm\'ediat que l'\'equation d'Euler-Lagrange pour ce probl\`eme variationnel
est exactement (\ref{cara8}). Tout consiste donc \`a savoir si ce probl\`eme
variationnel a une unique solution ou non. En g\'en\'eral, cela sera vrai pour
$A$ variant dans un domaine ouvert. Si tel est le cas, nous noterons

$${\cal Z} : (y,\epsilon,\pi)\hbox{ ou }(y,A) \longmapsto z,$$
l'application qui, \`a $(y,A)$, associe la solution $z$ de (\ref{cara8}).

\noindent {\bf V\'erification de (\ref{cara9}).}
Si $z = {\cal Z}(y, A)$, il vient, en utilisant la d\'efinition de
$H$ et (\ref{cara8}),

$$H^1_2(y,z) = {\partial F\over \partial z^j_1}(y,z)z^j_2
= A_{j,k}z^j_2z_2^k + A_{j,n+2}z^j_2
= A_{j,n+2}z^j_2,$$
qui co\"{\i}ncide avec la premi\`ere identit\'e de (\ref{cara9}). On
v\'erifie de
m\^eme l'identit\'e pour $H^2_1$.\\

\noindent {\bf V\'erification de (\ref{cara10}) et d\'etermination de $w$.}
Nous choisissons

\begin{equation}\label{cara11}
w = W(y,{\cal Z}(y,A),A).
\end{equation}
En utilisant (\ref{cara8}) et (\ref{cara11}), nous avons, notant $z = {\cal
Z}(y,A)$,

$$\begin{array}{ccl}
H^1_1(y,z) & = & \displaystyle {\partial F\over \partial z^j_1}(y,z)z^j_1 -
F(y,z)
= A_{j,k}z^j_1z_2^k + A_{j,n+2}z^j_1 - F(y,z)\\
& = & \left( A_{j,k}z_1^jz_2^k + A_{j,n+2}z_1^j
+A_{n+1,k}z_2^k + A_{n+1,n+2} - F(y,z)\right) 
- A_{n+1,k}z_2^k - A_{n+1,n+2}\\
& = & w - A_{n+1,k}z_2^k - A_{n+1,n+2}
\end{array}$$
et
$$\begin{array}{ccl}
H^2_2(y,z) & = & \displaystyle {\partial F\over \partial z^k_2}(y,z)z^k_2 -
F(y,z)
= A_{j,k}z_1^jz^k_2 + A_{n+2,k}z^k_2 - F(y,z)\\
& = & \left( A_{j,k}z_1^jz_2^k + A_{j,n+2}z_1^j
+A_{n+1,k}z_2^k + A_{n+1,n+2} - F(y,z)\right)
- A_{j,n+2}z_1^j - A_{n+1,n+2}\\
& = & w - A_{j,n+2}z_1^j - A_{n+1,n+2}.
\end{array}$$
Et nous obtenons ainsi (\ref{cara10}).

\subsubsection{Equations de Hamilton}
Supposant que l'\'equation
${\partial W\over \partial z}(y,z,\epsilon,\pi)=0$
a une unique solution $z = {\cal Z}(y,\epsilon,\pi)$, nous d\'efinissons l'hamiltonien ${\cal H}$ par

$${\cal H}: (y,\epsilon,\pi) \longmapsto W(y,{\cal Z}(y,\epsilon,\pi),\epsilon,\pi).$$
(${\cal H}(y,\epsilon,\pi)$ est aussi \'egal \`a $w$, solution de (\ref{cara10})).
Remarquons que nous avons les relations suivantes

\begin{equation}\label{cara12}
\begin{array}{cll}
\displaystyle {\partial {\cal H}\over \partial y^j}(y,\epsilon,\pi) & =
&\displaystyle 
{\partial W\over \partial y^j}(y,{\cal Z}(y,\epsilon,\pi),\epsilon,\pi) +
{\partial W\over \partial z_{\alpha}^k}(y,{\cal
Z}(y,\epsilon,\pi),\epsilon,\pi)
{\partial {\cal Z}_{\alpha}^k\over \partial y^j}(y,\epsilon,\pi)\\
& = & \displaystyle {\partial W\over \partial y^j}(y,{\cal
Z}(y,\epsilon,\pi),\epsilon,\pi)
= - {\partial F\over \partial y^j}(y,{\cal Z}(y,\epsilon,\pi)),
\end{array}
\end{equation}
et en notant
\begin{equation}\label{p1}
\left( \begin{array}{cc}
{\cal P}^1_1 & {\cal P}^1_2 \\ {\cal P}^2_1 & {\cal
P}^2_2\end{array}\right) (y,\epsilon,\pi)
= \left( \begin{array}{cc}
\pi^2_k{\cal Z}^k_2(y,\epsilon,\pi)+\epsilon_2^2 &
-(\pi^1_k{\cal Z}^k_2(y,\epsilon,\pi)+\epsilon_2^1)\\
-(\pi^2_j{\cal Z}^j_1(y,\epsilon,\pi)+\epsilon_1^2) &
\pi^1_j{\cal Z}^j_1(y,\epsilon,\pi)1+\epsilon_1^1
\end{array}\right) ,
\end{equation}

\noindent (remarquer qu'alors $W(y,z,\pi ,\epsilon ) =
{\cal P}^1_{\alpha}(y,\epsilon,\pi)\left( \pi ^{\alpha}_jz^j_1+\epsilon ^{\alpha}_1\right)
- F(y,z)
= {\cal P}^2_{\alpha}(y,\epsilon,\pi)\left( \pi ^{\alpha}_jz^j_2+\epsilon
^{\alpha}_2\right) -F(y,z)$)

\begin{equation}\label{cara13}
\begin{array}{cll}
\displaystyle {\partial {\cal H}\over \partial \pi
^{\alpha}_j}(y,\epsilon,\pi) & = &
\displaystyle {\partial W\over \partial \pi ^{\alpha}_j}
(y,{\cal Z}(y,\epsilon,\pi),\epsilon,\pi) +
{\partial W\over \partial z_{\beta}^k}(y,{\cal
Z}(y,\epsilon,\pi),\epsilon,\pi)
{\partial {\cal Z}_{\beta}^k\over \partial \pi ^{\alpha}_j}(y,\epsilon,\pi )\\
& = & \displaystyle {\partial W\over \partial \pi ^{\alpha}_j}
(y,{\cal Z}(y,\pi ,\epsilon ),\epsilon,\pi)
= {\cal P}^{\beta}_{\alpha}(y,\epsilon,\pi){\cal Z}^j_{\beta}(y,\epsilon,\pi),
\end{array}
\end{equation}

\begin{equation}\label{cara14}
\begin{array}{cll}
\displaystyle {\partial {\cal H}\over \partial \epsilon ^{\alpha}_{\beta}}
(y,\epsilon,\pi) & = &\displaystyle 
{\partial W\over \partial \epsilon ^{\alpha}_{\beta}}
(y,{\cal Z}(y,\epsilon,\pi),\epsilon,\pi) +
{\partial W\over \partial z_{\gamma}^k}(y,{\cal
Z}(y,\epsilon,\pi),\epsilon,\pi)
{\partial {\cal Z}_{\gamma}^k\over \partial \epsilon
^{\alpha}_{\beta}}(y,\epsilon,\pi)\\
& = & \displaystyle {\partial W\over \partial \epsilon
^{\alpha}_{\beta}}(y,{\cal Z}
(y,\epsilon,\pi),\epsilon,\pi)
= {\cal P}^{\beta}_{\alpha}(y,\epsilon,\pi).
\end{array}
\end{equation}

A pr\'esent, nous consid\'erons la fonctionnelle suivante.
A toute application $(u,\epsilon,\pi)$ d\'efinie sur un ouvert $\Omega$ de $\Bbb{C}$, nous associons
son graphe $\Sigma:=\{(t,u(t),\epsilon(t),\pi(t))/t\in \Omega\}$ et nous posons 

$$\begin{array}{cll}
{\cal A}_{\cal H}(u,\epsilon,\pi) & = & \displaystyle 
\int _{\Sigma} \left( \pi ^1_jdu^j+\epsilon ^1_{\alpha}dt^{\alpha}\right)
\wedge
\left( \pi ^2_kdu^k+\epsilon ^2_{\beta}dt^{\beta}\right) 
-{\cal H}(u,\epsilon,\pi)dt^1\wedge dt^2 \\
& = & \displaystyle \int _{\Omega}\left( \left| \begin{array}{cc}
\pi^1_j(t){\partial u^j\over \partial t^1}(t)+\epsilon_1^1(t) &
\pi^1_k(t){\partial u^k\over \partial t^2}(t)+\epsilon_2^1(t)\\
\pi^2_j(t){\partial u^j\over \partial t^1}(t)+\epsilon_1^2(t) &
\pi^2_k(t){\partial u^k\over \partial t^2}(t)+\epsilon_2^2(t)
\end{array}\right| -{\cal H}(u(t),\epsilon(t),\pi(t))\right) dt^1\wedge
dt^2.
\end{array}$$

Ecrivons les \'equations satisfaites par les points critiques de cette fonctionnelle. Pour cela, nous noterons
\begin{equation}\label{p2}
\left( \begin{array}{cc}
\hat{P}^1_1(t) & \hat{P}^1_2(t) \\ \hat{P}^2_1(t) & \hat{P}^2_2(t)\end{array}\right) =
\left( \begin{array}{cc}
\pi^2_k(t){\partial u^k\over \partial t^2}(t)+\epsilon_2^2(t) &
-(\pi^1_k(t){\partial u^k\over \partial t^2}(t)+\epsilon_2^1(t))\\
-(\pi^2_j(t){\partial u^j\over \partial t^1}(t)+\epsilon_1^2(t)) &
\pi^1_j(t){\partial u^j\over \partial t^1}(t)+\epsilon_1^1(t)
\end{array}\right) ,
\end{equation}

\noindent {\bf Variations par rapport \`a $u$}:

\begin{equation}\label{cara15}
{\partial \left( \hat{P}^{\alpha}_{\beta}\pi ^{\beta}_j\right) \over
\partial t^{\alpha}} = 
- {\partial {\cal H}\over \partial y^j}.
\end{equation}

\noindent {\bf Variations par rapport \`a $\pi$}: 

\begin{equation}\label{cara16}
\hat{P}^{\beta}_{\alpha}{\partial u^j\over \partial t^{\beta}} = 
{\partial {\cal H}\over \partial \pi ^{\alpha}_j}.
\end{equation}

\noindent {\bf Variations par rapport \`a $\epsilon$}:

\begin{equation}\label{cara17}
\hat{P}^{\beta}_{\alpha} = 
{\partial {\cal H}\over \partial \epsilon ^{\alpha}_{\beta}}.
\end{equation}

En comparant (\ref{cara13}) avec (\ref{cara16}) et (\ref{cara14}) avec
(\ref{cara17}), on obtient respectivement
${\cal P}^{\beta}_{\alpha}(y,\epsilon,\pi){\cal Z}^j_{\beta}(y,\epsilon,\pi) =
\hat{P}^{\beta}_{\alpha}{\partial u^j\over \partial t^{\beta}}$ et
${\cal P}^{\beta}_{\alpha}(y,\epsilon,\pi) = \hat{P}^{\beta}_{\alpha}$.
De ces deux \'equations, il vient ais\'ement

$${\partial u^k\over \partial t^{\alpha}} = {\cal Z}^k_{\alpha}(u(t),\epsilon(t),\pi(t)).$$
Cela entra\^{\i}ne que
$\hat{P}^{\alpha}_{\beta}\pi ^{\beta}_j = {\partial F\over \partial
z^j_{\alpha}}$.
En reportant cela dans (\ref{cara15}) et en utilisant (\ref{cara12}), on retrouve que $u$ est solution de
l'\'equation d'Euler-Lagrange (\ref{4.8}).

Nous pouvons aussi \'ecrire que l'action est stationnaire sous l'effet de
variations
de $t$, vue comme variable ind\'ependante (il faut penser ${\cal A}_{\cal H}$ comme
l'int\'egrale d'une 2-forme sur une surface plong\'ee dans l'espace des coordonn\'ees
$(t,y,\epsilon,\pi)$). Par exemple,

$$0 = \delta {\cal A}_{\cal H}(\delta t^1) =
\int _{\Sigma}\delta t^1 d\left( 
\epsilon ^2_1(\pi^1_jdu^j+\epsilon^1_{\alpha}dt^{\alpha})
- \epsilon ^1_1(\pi^2_kdu^k+\epsilon^2_{\beta}dt^{\beta})
+{\cal H}dt^2\right).$$
En calculant, on obtient

\begin{equation}\label{cara18}
{\partial \over \partial t^{\alpha}}\left( {\cal H}(u(t),\epsilon(t),\pi(t))\delta ^{\alpha}_{\beta}
- \hat{P}^{\alpha}_{\gamma}(t)\epsilon^{\gamma}_{\beta}(t)\right) =0.
\end{equation}
Cette \'equation exprime la conservation du tenseur \'energie-impulsion
${\partial H^{\alpha}_{\beta}\over \partial t^{\alpha}} = 0$, en vertu de
(\ref{cara7}).\\

Le syst\`eme d'\'equations (\ref{cara15}), (\ref{cara16}), (\ref{cara17})
et (\ref{cara18})
constitue un analogue des \'equations classiques de Hamilton, $(t,y)$
\'etant les
variables de temps et d'espace et $(\epsilon,\pi)$ celles d'\'energie et
d'impulsion.
Une forme diff\'erente est obtenue si, explo\^{\i}tant l'invariance par
transformation
de jauge, l'on choisit $(\epsilon, \pi)$ tels que

\begin{equation}\label{cara19}
{\partial \hat{P}^{\alpha}_{\beta}\over \partial t^{\alpha}} = 0.
\end{equation}
Pour voir cela, notons $\phi ^{\alpha} := \pi^{\alpha}_jdu^j +
\epsilon^{\alpha}_{\beta}dt^{\beta}$ et supposons que $\phi^1\wedge
\phi^2\neq 0$.
Alors, d'apr\`es un th\'eor\`eme de J. Moser [Mo] il existe des
coordonn\'ees sur $\Omega$
qui trivialisent la forme symplectique $\phi^1\wedge \phi^2$, c'est \`a
dire deux fonctions
$f^1,f^2:\Omega \longrightarrow \Bbb{R}$ telles que $\phi^1\wedge \phi^2 =
df^1\wedge df^2$.
Soit $g:\Omega \longrightarrow SL(2,\Bbb{R})$ l'unique fonction
telle que
$d\left( \begin{array}{c}f^1 \\ f^2\end{array}\right) =
g\left( \begin{array}{c}\phi^1 \\ \phi^2\end{array}\right)$. En rempla\c{c}ant
$(\epsilon, \pi)$ par $(\tilde{\epsilon}, \tilde{\pi}):=(g\epsilon, g\pi)$, on
obtient des variables d\'ecrivant le m\^eme probl\`eme lagrangien, mais
telles que 
$d\tilde{\phi}^1 = d\tilde{\phi}^2 = 0$, ce qui \'equivaut \`a (\ref{cara19}).

Dans un tel choix de coordonn\'ees , les \'equations (\ref{cara15}),
(\ref{cara16}), (\ref{cara17}) et (\ref{cara18}) donnent

\begin{equation}\label{cara20}
\begin{array}{cc}\displaystyle 
{\partial \left( \hat{P}^{\beta}_{\gamma}(t)\epsilon^{\gamma}_{\alpha}(t)-{\cal H}(u(t),\epsilon(t),\pi(t))\delta ^{\beta}_{\alpha}\right)
\over \partial t^{\beta}} = 0,
&\displaystyle 
{\partial \left( \hat{P}^{\beta}_{\alpha}(t)\pi ^{\alpha}_j(t)\right) \over
\partial t^{\beta}} = - {\partial {\cal H}\over \partial y^j}\\
\displaystyle 
{\partial \left( \hat{P}^{\beta}_{\alpha}(t)t^{\gamma}\right) \over
\partial t^{\beta}}= 
{\partial {\cal H}\over \partial \epsilon ^{\alpha}_{\gamma}},
&\displaystyle 
{\partial \left( \hat{P}^{\beta}_{\alpha}(t) u^j(t)\right) \over \partial
t^{\beta}} = 
{\partial {\cal H}\over \partial \pi ^{\alpha}_j}
\end{array}
\end{equation}

\subsubsection{Formulation avec contrainte}
Le param\`etre $w$, introduit initialement pour construire la correspondance de
Legendre-Carath\'eodory ne joue aucun r\^ole dans la dynamique. Nous pouvons imposer
que $w$ soit \'egal \`a une constante arbitraire $C$ et restreindre la correspondance \`a
$\{ (y,z,C)/(y,z)\in T{\cal N}^{\Bbb{C}}\}$. L'image de cette correspondance est alors
l'hypersurface

$${\cal H}^C:= \{ (y,\epsilon,\pi)/{\cal H}(y,\epsilon,\pi) = C\}.$$
Le probl\`eme variationnel est alors simplement obtenu en travaillant dans l'ensemble
des surfaces $\Sigma=\{ (t,u(t),\epsilon(t),\pi(t))/t\in \Omega\}$ plong\'ees dans
$\Omega \times {\cal H}^C$ avec la fonctionnelle

$$
{\cal A}(u,\epsilon,\pi) = 
\int _{\Sigma}
\left( \pi ^1_jdu^j+\epsilon ^1_{\alpha}dt^{\alpha}\right)
\wedge
\left( \pi ^2_kdu^k+\epsilon ^2_{\beta}dt^{\beta}\right) .$$
Les points critiques de ${\cal A}$ sous la contrainte
${\Sigma}\subset \Omega \times {\cal H}^C$
satisfont aux \'equations suivantes (avec la notation (\ref{p2}))

\begin{equation}\label{carac1}
\begin{array}{cc}\displaystyle 
{\partial \left( \hat{P}^{\beta}_{\gamma}(t)\epsilon^{\gamma}_{\alpha}(t)\right)
\over \partial t^{\beta}} = 0,
&\displaystyle 
{\partial \left( \hat{P}^{\beta}_{\alpha}(t)\pi ^{\alpha}_j(t)\right) \over
\partial t^{\beta}} = - \mu{\partial {\cal H}\over \partial y^j}\\
\displaystyle 
\hat{P}^{\beta}_{\alpha}(t){\partial t^{\gamma}\over
\partial t^{\beta}}= 
\mu{\partial {\cal H}\over \partial \epsilon ^{\alpha}_{\gamma}},
&\displaystyle 
\hat{P}^{\beta}_{\alpha}(t){\partial u^j(t)\over \partial
t^{\beta}} = 
\mu{\partial {\cal H}\over \partial \pi ^{\alpha}_j},
\end{array}
\end{equation}
o\`u $\mu$ est le multiplicateur de Lagrange associ\'e \`a la contrainte. En utilisant
(\ref{cara13}) et (\ref{cara14}), on en d\'eduit

\begin{equation}\label{carac2}
\mu{\cal P}^{\beta}_{\alpha}(u,\epsilon,\pi){\cal Z}^j_{\beta}(u,\epsilon,\pi) =
\hat{P}^{\beta}_{\alpha}{\partial u^j\over \partial t^{\beta}}
\end{equation}
et
\begin{equation}\label{carac3}
\mu{\cal P}^{\beta}_{\alpha}(u,\epsilon,\pi) = \hat{P}^{\beta}_{\alpha}.
\end{equation}
En substituant (\ref{carac3}) dans (\ref{carac2}), on obtient
$\hat{P}^{\beta}_{\alpha}{\partial u^j\over \partial t^{\beta}} =
\hat{P}^{\beta}_{\alpha}{\cal Z}^j_{\beta}(u,\epsilon,\pi)$, d'o\`u l'on tire
${\partial u^j\over \partial t^{\beta}} ={\cal Z}^j_{\beta}(u,\epsilon,\pi)$. En reportant dans
(\ref{carac3}), on conclut que $\mu=1$. Le syst\`eme (\ref{carac1}) conduit donc \`a

\begin{equation}\label{carac4}
\begin{array}{cc}\displaystyle 
{\partial \left( \hat{P}^{\beta}_{\gamma}(t)\epsilon^{\gamma}_{\alpha}(t)\right)
\over \partial t^{\beta}} = 0,
&\displaystyle 
{\partial \left( \hat{P}^{\beta}_{\alpha}(t)\pi ^{\alpha}_j(t)\right) \over
\partial t^{\beta}} = - {\partial {\cal H}\over \partial y^j}\\
\displaystyle 
\hat{P}^{\beta}_{\alpha}(t){\partial t^{\gamma}\over
\partial t^{\beta}}= 
{\partial {\cal H}\over \partial \epsilon ^{\alpha}_{\gamma}},
&\displaystyle 
\hat{P}^{\beta}_{\alpha}(t){\partial u^j(t)\over \partial
t^{\beta}} = 
{\partial {\cal H}\over \partial \pi ^{\alpha}_j}.
\end{array}
\end{equation}

\subsubsection{Propri\'et\'es du hamiltonien pour un probl\`eme invariant par 
transformation conforme}

Nous notons $\vec{A}:=\left( A_{j,k}\right) _{1\leq j,k\leq n}$. Il est
possible de
r\'e\'ecrire la fonctionnelle $W$ sous la forme

$$
W(y,z,A) = -{i\over 2}\vec{A}_{jk}\overline{z}^jz^k
+ A_{j,n+2}z_1^j + A_{n+1,k}z_2^k +A_{n+1,n+2} - F(y,z).
$$

Donc, les \'equations (\ref{cara8}), \`a r\'esoudre pour trouver $z$ en
fonction
de $y$ et de $A$, sont \'equivalentes \`a la relation suivante

\begin{equation}\label{cara21}
{\partial W\over \partial \overline{z}^j} =
- {\partial F\over \partial \overline{z}^j}
-{i\over 2} \vec{A}_{jk}z^k +
{1\over 2} \left(  A_{j,n+2} + iA_{n+1,j}\right) =0.
\end{equation}
En multipliant cette \'equation par $\overline{z}^j$, en sommant et en utilisant
(\ref{4.3}),
nous obtenons 

$$\begin{array}{rcc}
\left[ - F(y,z) -{i\over 2} \vec{A}_{jk}\overline{z}^jz^k +
{1\over 2} \left(  A_{j,n+2}z_1^j + A_{n+1,j}z_2^j\right) \right] 
+{i\over 2} \left[ A_{n+1,j}z_1^j - A_{j,n+2}z_2^j\right] & = & \\
\displaystyle - {\partial F\over \partial \overline{z}^j}\overline{z}^j 
-{i\over 2} \vec{A}_{jk}\overline{z}^jz^k +
{1\over 2} \left(  A_{j,n+2} + iA_{n+1,j}\right) \overline{z}^j & = &
0.
\end{array}$$
On en conclut que

\begin{equation}\label{cara22a}
A_{n+1,j}z_1^j - A_{j,n+2}z_2^j = 0
\end{equation}
et en reportant dans l'expression pr\'ec\'edente de $W(y,z,A)$,

\begin{equation}\label{cara22}
{\cal H}(y,A) = {1\over 2} \left(  A_{j,n+2}z_1^j + A_{n+1,j}z_2^j\right) + A_{n+1,n+2}.
\end{equation}

Nous pouvons d\'eduire du Lemme 1 et de (\ref{cara21}) que, $\vec{A}_{jk}$
\'etant
donn\'e, $z$ est une fonction homog\`ene complexe de degr\'e 1 de
$A_{j,n+2} + iA_{n+1,j}$,
id est (notant ${\cal Z} = {\cal Z}_1+i{\cal Z}_2$)

$${\cal Z}(y,\vec{A},\lambda(A_{j,n+2}+iA_{n+1,j})) =
\lambda {\cal Z}(y,\vec{A},A_{j,n+2}+iA_{n+1,j}),$$
pour tout $\lambda \in \Bbb{C}$; (\ref{cara22}) entra\^{\i}ne alors que
${\cal H}(y,z,A) = \tilde{\cal H}(y,\vec{A},A_{j,n+2}+iA_{n+1,j}) +
A_{n+1,n+2}$, o\`u $\tilde{\cal H}$ est tel que

$$\tilde{\cal H}(y,\vec{A},\lambda (A_{j,n+2}+iA_{n+1,j})) =
|\lambda |^2\tilde{\cal H}(y,\vec{A},A_{j,n+2}+iA_{n+1,j}).$$

Dans les variables $(\epsilon,\pi)$, cela signifie que

$${\cal H}(y,\lambda(\epsilon^{\alpha}_1+i\epsilon^{\alpha}_2),\pi) =
|\lambda |^2{\cal H}(y,\epsilon^{\alpha}_1+i\epsilon^{\alpha}_2,\pi).$$

\subsubsection{Le cas hermitien}
A titre d'exemple, nous examinons ce qui se passe lorsque
$F(y,z) = {1\over 2}h_{jk}(y)\overline{z^j}z^k$, o\`u
$h_{jk}(y)=g_{jk}(y)-i\omega_{jk}(y)$
est un tenseur m\'etrique hermitien. Etudions d'abord la correspondance de
Legendre-Carath\'eodory.
Nous avons alors l'expression suivante pour $W$

$$
W(y,z,A) = 
-{1\over 2}\left[ g_{jk} -i\left( \omega_{jk} - \vec{A}_{jk}\right) \right]
\overline{z}^jz^k
+ A_{j,n+2}z_1^j + A_{n+1,k}z_2^k +A_{n+1,n+2}.
$$

Alors les \'equations (\ref{cara8}) sont \'equivalentes \`a

$${\partial W\over \partial \overline{z}^j} =
-{1\over 2} \left( h_{jk} +i \vec{A}_{jk}\right) z^k +
{1\over 2} \left(  A_{j,n+2} + iA_{n+1,j}\right) =0.$$

Et le probl\`eme est donc de trouver $z$, solution de

\begin{equation}\label{cara23}
h^{\star}_{jk}z^k = A_{j,n+2} + iA_{n+1,j},
\end{equation}
avec $h^{\star}_{jk}:= h_{jk} +i \vec{A}_{jk}$. Cette \'equation a une
unique solution
si et seulement si det$h^{\star}_{jk} \neq 0$.

\begin{rema}
Il est possible d'\'ecrire tout le syst\`eme (\ref{cara8}), (\ref{cara9}),
(\ref{cara10})
sous une forme condens\'ee, de la mani\`ere suivante. Nous notons

$$\zeta_1:= \left( \begin{array}{c}z_1^1\\ \vdots
\\z_1^n\\1\\0\end{array}\right) \in \Bbb{R}^{n+2},
\zeta_2:= \left( \begin{array}{c}z_2^1\\ \vdots
\\z_2^n\\0\\1\end{array}\right) \in \Bbb{R}^{n+2},
\zeta := \zeta_1+i\zeta_2 \in \Bbb{C}^{n+2}.$$
Alors $W(y,z,A) = -{i\over 2}A_{jk}\overline{\zeta}^j\zeta^k - F(y,z)$ et
(\ref{cara8}), (\ref{cara9}), (\ref{cara10}) est \'equivalent \`a

$$G\zeta = 0,$$
o\`u 

$$G:= \left( \begin{array}{cc}h & 0\\0 & w\1 _2-{^tH}\end{array}\right) +iA
= \left( \begin{array}{ccccc}
h^{\star}_{11} & \dots & h^{\star}_{1n} & iA_{1,n+1} & iA_{1,n+2}\\
\vdots & & \vdots & \vdots & \vdots \\
h^{\star}_{n1} & \dots & h^{\star}_{nn} & iA_{n,n+1} & iA_{n,n+2}\\
iA_{n+1,1} & \dots & iA_{n+1,n} & w -H^1_1 & -H^2_1+iA_{n+1,n+2}\\
iA_{n+2,1} & \dots & iA_{n+2,n} & -H^1_2+iA_{n+2,n+1} & w -H^2_2\\
\end{array}\right) .$$
Ici, $h_{jk}$ et $A$ sont donn\'es et $z$, $w$ et $H$ sont les inconnues.
\end{rema}

Dans cette situation, nous pouvons expliciter d'avantage ${\cal H}$: notant
$K_{\vec{A}}$ l'inverse de $h^{\star} = h+i\vec{A}$, c'est \`a dire tel que

$$K_{\vec{A}}^{jk}\left( h_{kl}+i\vec{A}_{kl}\right) = \delta^j_l,$$
nous d\'eduisons de (\ref{cara23})

$$z^j = K_{\vec{A}}^{jk}\left( A_{k,n+2}+iA_{n+1,k}\right) ,$$
et en reportant dans (\ref{cara22}) (explo\^{\i}tant (\ref{cara22a})),

\begin{equation}\label{cara24}
{\cal H}(y,A) = {1\over 2}K_{\vec{A}}^{jk}\left( A_{j,n+2}-iA_{n+1,j}\right) 
\left( A_{k,n+2}+iA_{n+1,k}\right) + A_{n+1,n+2}.
\end{equation}

\subsubsection{Une g\'en\'eralisation}
On peut introduire, \`a la place des variables $\epsilon$, $\pi$, une famille
de variables $\epsilon^{(J)}=\epsilon^{(J)\alpha}_{\beta}$ et
$\pi^{(J)}=\pi^{(J)\alpha}_j$, o\`u $J=1,...,N$, pour un certain entier $N$.
On remplace alors la d\'efinition pr\'ec\'edente de $W$ par

$$W(y,z,\epsilon^{(J)},\pi^{(J)}) := \sum _{J=1}^N
\left| \begin{array}{cc}
\pi^{(J)1}_jz^j_1+\epsilon_1^{(J)1} & \pi^{(J)1}_kz^k_2+\epsilon_2^{(J)1}\\
\pi^{(J)2}_jz^j_1+\epsilon_1^{(J)2} & \pi^{(J)2}_kz^k_2+\epsilon_2^{(J)2}
\end{array}\right| - F(y,z).
$$
A $y$ fix\'e, les variables $(z,w)$ sont li\'ees aux variables
$(\epsilon^{(J)},\pi^{(J)})$
de fa\c{c}on telle que $z$ soit la solution ${\cal Z}(y,\epsilon^{(J)},\pi^{(J)})$ de
${\partial W\over \partial z}(y,z,\epsilon^{(J)},\pi^{(J)})=0$ et
$w = W(y,{\cal Z}(y,\epsilon^{(J)},\pi^{(J)}),\epsilon^{(J)},\pi^{(J)})$.
L'analyse de cette correspondance de Legendre-Carath\'eodory est identique \`a ce qui
pr\'ec\`ede, il
suffit de prendre comme nouvelle d\'efinition de $\left( A_{jk}\right)
_{1\leq j,k\leq n+2}$:

$$\begin{array}{c}
\displaystyle
A_{j,k}:= \sum _{J=1}^N
\left| \begin{array}{cc}
\pi^{(J)1}_j & \pi^{(J)1}_k\\
\pi^{(J)2}_j & \pi^{(J)2}_k
\end{array}\right| ,\
A_{j,n+2}:=\sum _{J=1}^N\left| \begin{array}{cc}
\pi^{(J)1}_j & \epsilon_2^{(J)1}\\
\pi^{(J)2}_j & \epsilon_2^{(J)2}
\end{array}\right| ,\\
\displaystyle
A_{n+1,k}:=\sum _{J=1}^N\left| \begin{array}{cc}
\epsilon_1^{(J)1} & \pi^{(J)1}_k\\
\epsilon_1^{(J)2} & \pi^{(J)2}_k
\end{array}\right| \hbox{ et }
A_{n+1,n+2}:=\sum _{J=1}^N\left| \begin{array}{cc}
\epsilon_1^{(J)1}& \epsilon_2^{(J)1}\\
\epsilon_1^{(J)2}& \epsilon_2^{(J)2}
\end{array}\right| .
\end{array}$$
Il ressort de l'analyse faite que g\'en\'eriquement, \`a tout
$(\epsilon^{(J)},\pi^{(J)})$,
correspond un unique $(z,w)$. De plus, ${\cal Z}(y,\epsilon^{(J)},\pi^{(J)})$ est maintenant
une fonction invariante par le groupe des transformations symplectiques de
$\Bbb{R}^{2N}$,
pr\'eservant la 2-forme $dx^{(1)1}\wedge dx^{(1)2} +...+dx^{(N)1}\wedge
dx^{(N)2}$. On d\'efinit ainsi
${\cal H}(y,\epsilon^{(J)},\pi^{(J)}) =
W(y,{\cal Z}(y,\epsilon^{(J)},\pi^{(J)}),\epsilon^{(J)},\pi^{(J)})$.

Les solutions du probl\`eme variationnel initial peuvent \^etre
obtenues comme suit. Nous associons \`a chaque application $(u,\epsilon^{(J)},\pi^{(J)})$
son graphe $\Sigma:= \{(t,u(t),\epsilon^{(J)}(t),\pi^{(J)}(t))/t\in \Omega\}$ et nous
utilisons l'un des deux probl\`emes variationnels suivants

$${\cal A}_{\cal H}(u,\epsilon^{(J)},\pi^{(J)}) =  
\int _{\Sigma}
\sum _{J=1}^N\left( \pi ^{(J)1}_jdu^j
+\epsilon ^{(J)1}_{\alpha}dt^{\alpha}\right) \wedge
\left( \pi ^{(J)2}_kdu^k+\epsilon ^{(J)2}_{\beta}dt^{\beta}\right) 
-{\cal H}(u,\pi^{(J)},\epsilon^{(J)})dt^1\wedge dt^2;$$
ou bien 

$${\cal A}(u,\epsilon^{(J)},\pi^{(J)}) =  
\int _{\Sigma}
\sum _{J=1}^N\left( \pi ^{(J)1}_jdu^j
+\epsilon ^{(J)1}_{\alpha}dt^{\alpha}\right) \wedge
\left( \pi ^{(J)2}_kdu^k+\epsilon ^{(J)2}_{\beta}dt^{\beta}\right) ,$$
avec la contrainte ${\cal H}(y,\epsilon^{(J)},\pi^{(J)})=C$.
Notant

$$\left( \begin{array}{cc}
\hat{P}^{(J)1}_1 & \hat{P}^{(J)1}_2 \\ \hat{P}^{(J)2}_1 &
\hat{P}^{(J)2}_2\end{array}\right) =
\left( \begin{array}{cc}
\pi^{(J)2}_k{\partial u^k\over \partial t^2}+\epsilon_2^{(J)2} &
-(\pi^{(J)1}_k{\partial u^k\over \partial t^2}+\epsilon_2^{(J)1})\\
-(\pi^{(J)2}_j{\partial u^j\over \partial t^1}+\epsilon_1^{(J)2}) &
\pi^{(J)1}_j{\partial u^j\over \partial t^1}+\epsilon_1^{(J)1}
\end{array}\right) ,$$
on obtient, pour les points critique de ${\cal A}$ sous la contrainte
${\cal H}(y,\epsilon^{(J)},\pi^{(J)})=C$, le syst\`eme d'\'equations 
d'Euler-Lagrange suivant

\begin{equation}\label{cara25}
\begin{array}{cc}\displaystyle 
{\partial \over \partial t^{\beta}}\left( \sum _{J=1}^N
\hat{P}^{(J)\beta}_{\gamma}(t)\epsilon^{\gamma}_{\alpha}(t)
\right) = 0,
&\displaystyle 
{\partial \over \partial t^{\beta}}\left( \sum _{J=1}^N\hat{P}^{(J)\beta}_{\alpha}(t)\pi
^{\alpha}_j(t)\right) = - {\partial {\cal H}\over \partial y^j}\\
\displaystyle 
\hat{P}^{(J)\beta}_{\alpha}(t) = 
{\partial {\cal H}\over \partial \epsilon ^{(J)\alpha}_{\beta}},
&\displaystyle 
\hat{P}^{(J)\beta}_{\alpha}(t) {\partial u^j\over \partial t^{\beta}}(t) = 
{\partial {\cal H}\over \partial \pi ^{(J)\alpha}_j},
\end{array}
\end{equation}
analogue aux \'equations de Hamilton.

Il semble int\'eressant d'essayer de comprendre les propri\'et\'es de ce
genre de probl\`eme.
On peut observer que dans le cas o\`u $N=1$, qui correspond au cas o\`u $A$
est de rang
deux, \`a une valeur $(y,z,w)$ est associ\'ee, via la correspondance de Legendre-Carath\'edodory, 
une famille $\{ (y,g\epsilon,g\pi)/g\in SL(2,\Bbb{R})$, mais dans le cas
o\`u $N>1$,
$A$ peut \^etre de n'importe quel rang, compris entre 2 et $2N$ et, une
m\^eme valeur $(y,z,w)$ correspond \`a plusieurs familles de valeurs de
$(\epsilon^{(J)},\pi^{(J)})$, selon le
rang de $A$. Cela rappelle une situation bien connue en m\'ecanique
quantique: l'espace des \'etats en m\'ecanique quantique est plus gros
que l'espace des \'etats de la m\'ecanique classique et  autorise la superpositions d'\'etats quantiques purs.

\subsection{Equations d'Hamilton-Jacobi}
Historiquement, la d\'emarche de Carath\'eodory fut de construire une g\'en\'eralisation
de l'\'equation d'Hamilton-Jacobi pour les probl\`emes variationnels \`a plusieurs
variables, ce qui l'a amen\'e assez naturellement \`a
d\'efinir l'hamiltonien ${\cal H}$. Plus tard H. Rund a \'ecrit un syst\`eme
d'\'equations canoniques associ\'ees \` a cet hamiltonien. C'est pourquoi, il me
semble int\'eressant de rappeler les \'equations d'Hamilton-Jacobi pour les diff\'erentes
th\'eories que nous avons rencontr\'ees.\\

Consid\'erons d'abord un probl\`eme variationnel \`a une variable. Soit
$I$ un intervalle de $\Bbb{R}$, $U$ un ouvert de $\Bbb{R}^n$ et
$L:I\times U\times \Bbb{R}^n\longrightarrow \Bbb{R}$ un lagrangien satisfaisant
la condition de Legendre. Nous cherchons une fonction $S:I\times U\longrightarrow \Bbb{R}$,
telle que $\forall (t,y,z)\in I\times U\times \Bbb{R}^n$,

\begin{equation}\label{jaco1}
L(t,y,z)\geq {\partial S\over \partial y^j}(t,y)z^j
+ {\partial S\over \partial t}(t,y),
\end{equation}
avec \'egalit\'e pour un certain $z=\psi(t,y)$:

\begin{equation}\label{jaco2}
L(t,y,\psi(t,y))= {\partial S\over \partial y^j}(t,y)\psi^j(t,y)
+ {\partial S\over \partial t}(t,y).
\end{equation}
$\psi$ est appel\'e {\em champ de Mayer} \footnote{noter que pour toute application
$u:I\longrightarrow U$,
${\partial S\over \partial y^j}(t,u(t)){\partial u^j\over \partial t}(t)
+ {\partial S\over \partial t}(t,u(t)) = {d\over dt}\left( S(t,u(t))\right)$}.
En particulier, pour tout $(t,y)\in I\times U$ fix\'e, $\psi(t,y)$ est un
minimum de $z\longmapsto L(t,y,z) - {\partial S\over \partial y^j}(t,y)z^j
+ {\partial S\over \partial t}(t,y)$, ce qui entra\^\i ne

\begin{equation}\label{jaco3}
{\partial L\over \partial z^j}(t,y,\psi(t,y)) = {\partial S\over \partial y^j}(t,y).
\end{equation}
Cette relation a la cons\'equence que

\begin{equation}\label{jaco4}
H\left( t,y,{\partial S\over \partial y^j}(t,y)\right) =
{\partial L\over \partial z^j}(t,y,\psi(t,y))\psi^j(t,y) -  L(t,y,\psi(t,y)).
\end{equation}
Maintenant, en substituant (\ref{jaco3}) dans (\ref{jaco2}), on obtient

$$L(t,y,\psi(t,y))={\partial L\over \partial z^j}(t,y,\psi(t,y))\psi^j(t,y)
+ {\partial S\over \partial t}(t,y).$$
Cette derni\`ere relation signifie exactement, gr\^ace \`a (\ref{jaco4}) que
$S$ est solution de

\begin{equation}\label{jaco5}
H\left( t,y,{\partial S\over \partial y^j}(t,y)\right)
+ {\partial S\over \partial t}(t,y) =0,
\end{equation}
l'\'equation d'Hamilton-Jacobi.\\

Une \'equation analogue pour la th\'eorie de Weyl s'obtient comme suit.
Nous partons de $L:\Omega\times U\times \Bbb{R}^{mn}$ ($\Omega\subset \Bbb{R}^m$,
$U\subset \Bbb{R}^n$) satisfaisant la condition de Legendre et nous cherchons
$S:\Omega\times U\longrightarrow \Bbb{R}^n$ tel que,
$\forall (t,y,z)\in \Omega\times U\times \Bbb{R}^{mn}$,

\begin{equation}\label{jaco6}
L(t,y,z)\geq {\partial S^{\alpha}\over \partial y^j}(t,y)z^j_{\alpha}
+ {\partial S^{\alpha}\over \partial t^{\alpha}}(t,y),
\end{equation}
avec \'egalit\'e pour un certain $z=\psi(t,y)$. Par le m\^eme raisonnement, on
trouve comme condition n\'ecessaire et suffisante sur $S$:

\begin{equation}\label{jaco7}
H\left( t,y,{\partial S^{\alpha}\over \partial y^j}(t,y)\right)
+ {\partial S^{\alpha}\over \partial t^{\alpha}}(t,y) =0.
\end{equation}

Enfin, la th\'eorie de Carath\'edory correspond au probl\`eme suivant: trouver
$S:\Omega\times U\longrightarrow \Bbb{R}^n$ tel que,
$\forall (t,y,z)\in \Omega\times U\times \Bbb{R}^{mn}$,

\begin{equation}\label{jaco8}
L(t,y,z)\geq \det\left( {\partial S^{\alpha}\over \partial y^j}(t,y)z^j_{\beta}
+ {\partial S^{\alpha}\over \partial t^{\beta}}(t,y)\right) ,
\end{equation}
avec \'egalit\'e pour $z=\psi(t,y)$. On constate alors que
$\pi^{\alpha}_j={\partial S^{\alpha}\over \partial y^j}(t,y)$,
$\epsilon^{\alpha}_{\beta}={\partial S^{\alpha}\over \partial t^{\beta}}(t,y)$
forment une solution de (\ref{cara5}), (\ref{cara6}) et (\ref{cara7}) avec $z=\psi(t,y)$
et $w=0$. Cette derni\`ere condition $w=0$ se transcrit en l'\'equation
d'Hamilton-Jacobi

\begin{equation}\label{jaco9}
{\cal H}\left( t,y,{\partial S^{\alpha}\over \partial t^{\beta}}(t,y),
{\partial S^{\alpha}\over \partial y^j}(t,y)\right) =0.
\end{equation}
On peut envisager une g\'en\'eralisation sur des familles de $N$ fonctions
$S^{(J)}:\Omega\times U\longrightarrow \Bbb{R}^n$, pour $J=1,...,N$,
\`a partir du formalisme propos\'e en 5.2.6.

\section{R\'ef\'erences}

\begin{enumerate}
\item[[BCh1]] D. Bao, S.S. Chern, {\em On a notable connection in Finsler
geometry}, Houston J. Math. 19 (1993), 135-180.
\item[[BCh2]] D. Bao, S.S. Chern, {\em A note on the Gauss-Bonnet theorem
for Finsler spaces}, Ann. Math. 143 (1996), 233-252.
\item[[BChS]] D. Bao, S.S. Chern, Z. Shen, {\em Finsler geometry},
(proceedings of the joint summer research conference om Finsler geometry)
Cont. Math., vol. 196, Amer. Math. Soc..
\item[[Br]] R. Bryant, {\em Finsler surfaces with prescribed curvature
conditions}, pr\'epublication, \`a para\^\i tre dans les Aisenstadt lectures
de R. Bryant, voir\\
{\em http://www.math.duke.edu/faculty/bryant/}
\item[[Ca]] C. Carath\'eodory, {\em Calculus of variations and partial differential
equations of first order}, Parts I, II, Holden-Day, San Francisco (1967), traduction
de {\em Variationsrechnung und partielle differentialgleichungen erster Ordnung},
Teubner, Leipzig und Berlin (1935).
\item[[Cn]] E. Cartan, {\em Les espaces m\'etriques fond\'es sur la notion
d'aire}, Actualit\'es scientifiques 72, Paris (1933).
\item[[Ch]] S.S. Chern, {\em Finsler geometry is just Riemannian geometry
without the quadratic restriction}, Notices of the AMS, September 1996,
959-963.
\item[[Fo]] P. Foulon, {\em G\'eom\'etrie des \'equation diff\'erentielles
du second ordre}, Ann. Inst. Henri Poincar\'e 45 (1) (1986), 1-28.
\item[[GiHi]] M. Giaquinta, S. Hildebrandt, {\em Calculus of variations II},
Grundlehren der mathematischen Wissenschaft 311, Springer 1996.
\item[[H\'e]] F. H\'elein, {\em Applications harmoniques, lois de 
conservation et rep\`eres mobiles}, Diderot \'editeur, Paris 1996; or {\em 
Harmonic maps, conservation laws and moving frames}, Diderot \'editeur, 
Paris 1997.
\item[[HK]] F. H\'elein, J. Kouneiher, {\em Hamiltonian formalism with several variables and quantum field theory, Part I},
math-ph/0004020 et {\em Finite dimensional Hamiltonian formalism for gauge and field theories},
math-ph/0010036.
\item[[Ma]] D. H. Martin, {\em Canonical variables and geodesic fields for the calculus of
variations of multiple integrals in parametric form}, Math. Z. 104 (1968), 16-27.
\item[[Mo]] J. Moser, {\em On the volume elements of a manifold}, Trans. Am. Math. Soc. 120 (1956), 286-294.
\item[[Ru]] H. Rund, {\em The Hamilton-Jacobi theory in the calculus of
variations, its role in Mathematics and
Physics}, Krieger Pub. 1973 (nouvelle \'edition avec un appendice
suppl\'ementaire).
\item[[Ru1]] H. Rund, {\em A canonical formalism for multiple integral problems in the
calculus of variations}, Aeq. Math. 3 (1968), 44-63.

\end{enumerate}

\end{document}